\documentclass[11pt,a4paper]{article}
\usepackage[utf8]{inputenc}
\usepackage{lmodern}
\usepackage{amsmath,amssymb,amsthm,mathtools}
\usepackage{booktabs,array,longtable}
\usepackage{microtype}
\usepackage[left=22mm,right=22mm,top=30mm,bottom=30mm]{geometry}
\usepackage{appendix}

\usepackage{enumitem,needspace}
\usepackage{xcolor,fancyhdr}
\usepackage{xurl}
\usepackage{hyperref}
\usepackage{bookmark}
\definecolor{linkblue}{RGB}{25,62,100}
\hypersetup{colorlinks=true,linkcolor=linkblue,
citecolor=linkblue, urlcolor=linkblue}
\usepackage{indentfirst}
\numberwithin{equation}{section}
\newtheorem{theorem}{Theorem}[section]
\newtheorem{proposition}[theorem]{Proposition}
\newtheorem{lemma}[theorem]{Lemma}
\newtheorem{conjecture}[theorem]{Conjecture}

\theoremstyle{definition}

\theoremstyle{remark}
\newtheorem{remark}[theorem]{Remark}
\newcommand{\Irr}{\operatorname{Irr}}
\newcommand{\IBr}{\operatorname{IBr}}
\newcommand{\Aut}{\operatorname{Aut}}
\newcommand{\Out}{\operatorname{Out}}
\newcommand{\Gal}{\operatorname{Gal}}
\newcommand{\GL}{\operatorname{GL}}
\newcommand{\SL}{\operatorname{SL}}
\DeclareMathOperator{\PSL}{PSL}
\DeclareMathOperator{\PSU}{PSU}

\DeclareMathOperator{\Inn}{Inn}
\DeclareMathOperator{\St}{St}

\newcommand{\Syl}{\operatorname{Syl}}

\newcommand{\dz}{\operatorname{dz}}
\newcommand{\Sym}{\operatorname{Sym}}
\newcommand{\Ind}{\operatorname{Ind}}
\newcommand{\Inf}{\operatorname{Inf}}
\newcommand{\Res}{\operatorname{Res}}

\newcommand{\OO}{\mathrm{O}}
\newcommand{\ZZ}{\mathrm{Z}}
\newcommand{\HH}{\mathcal H}
\newcommand{\WW}{\mathcal W}
\newcommand{\FF}{\mathbb F}
\newcommand{\QQ}{\mathbb Q}
\newcommand{\NN}{\mathrm N}
\newcommand{\CC}{\mathrm C}
\newcommand{\et}{\widetilde E}
\newcommand{\AT}{\widetilde A}

\newcommand{\mcite}[2]{\cite[#2]{#1}}

\title{\bfseries The inductive Navarro-Alperin weight condition\\ for the Suzuki and Ree groups}
\author{Qianhu Zhou}
\date{\scriptsize School of Mathematics and Statistics, Central China Normal University, Wuhan, 430079, China}
\begin{document}
\maketitle
{\noindent\small{\bf Abstract}
In this paper, we give verifications of  the inductive Navarro-Alperin weight condition in \cite{FFZ} for the Suzuki and Ree groups at all primes.}


\section{Introduction}\label{sec:intro}

In the representation theory of finite groups, the local-global conjectures postulate connections between the representation theory of a finite group and that of its local subgroups. Fix a prime $p$, and let $G$ be a finite group. A $p$-weight of $G$ is a pair $(R,\varphi)$, where $R$ is a $p$-subgroup of $G$ and $\varphi$ is an irreducible complex character of $\mathrm{N}_G(R)/R$ satisfying $\varphi(1)_p = |\mathrm{N}_G(R)/R|_p$. The Alperin weight conjecture \cite{Alperin1987} states that the number of irreducible $p$-Brauer characters of $G$ equals the number of $G$-conjugacy classes of $p$-weights of $G$.
In 2004, Navarro \cite{Navarro2004} refined this conjecture and predicted the existence of a bijection between the set of irreducible $p$-Brauer characters of $G$ and the set of $G$-conjugacy classes of $p$-weights of $G$ that may be chosen equivariant under certain Galois automorphisms.
In \cite{FFZ}, this refinement is termed the \textit{Navarro-Alperin weight (NAW) conjecture}, also known as the \textit{Galois Alperin weight conjecture} (see Conjecture \ref{con2.1} below). The present paper is devoted to this conjecture.

Feng, Fu and Zhou \cite{FFZ} work with a formulation of the NAW conjecture incorporating both Galois automorphisms and group automorphisms, and reduce the conjecture to a problem on simple groups. This gives rise to the so-called \textit{inductive Navarro-Alperin weight (NAW) condition} (see \cite[Definition 7.1]{FFZ}), which is typically stronger than the original conjecture. The NAW conjecture has been verified in several special cases, including $p$-solvable groups and finite groups with abelian Sylow $2$-subgroups \cite{FFZ}. Moreover, the inductive NAW condition has been established for all finite nonabelian simple groups of Lie type in their defining characteristic in \cite{FFZ}, as well as for the Ree simple groups ${}^2G_2(3^{2m+1})$ with $m\ge 1$ at all primes, and so on.

The main result of this paper is the following theorem.
\begin{theorem}\label{thm:main}
The inductive Navarro-Alperin weight condition is satisfied for  the Suzuki and Ree groups at all primes $p$.
\end{theorem}

The structure of this paper is as follows. 
In Section \ref{section 2}, we introduce some basic notations and a criterion for the inductive NAW condition of certain  non-abelian simple groups. In Section \ref{section 3}, we will verify that Suzuki simple groups and Ree simple groups with cyclic Sylow subgroups satisfy the inductive NAW condition. In Section \ref{section 4}, we treat  the case of the Ree groups with  non-cyclic Sylow subgroups.
In Section \ref{section 5}, we address several exceptional cases mentioned at the start of Subsection \ref{section 2.3}.
In section \ref{section 6}, we prove Theorem \ref{thm:main}. Finally, in Appendix \ref{app:maps}, we record
some computational steps needed for the exceptional situations encountered in Section \ref{section 5}.

\section{Preliminaries}\label{section 2}
We introduce the notation we use throughout this paper. 
\subsection{General notation}
Let $G$ be a finite group. Concerning the Brauer and ordinary character theory of $G$ we mainly follow the notation of \cite{N1998} and \cite{I1976}. For induction, restriction and inflation of characters we use the notation $\Ind$, $\Res$ 
and $\Inf$. Let $p$ be a prime. Write $\IBr_p(G)$ and $\Irr(G)$ for the sets of all irreducible $p$-Brauer characters and ordinary characters of $G$, respectively. If $B$ is a $p$-block of $G$, then write $\IBr_p(B)$ for the set of irreducible $p$-Brauer characters $\eta\in\IBr_p(G)$ such that $\eta$ belongs to the block $B$. A similar statement holds for ordinary characters. An irreducible ordinary character $\chi\in\Irr(G)$ is said to be \textit{$p$-defect-zero} if $\chi(1)_p=|G|_p$. Write $\dz_p(G)$ for the set of $p$-defect-zero characters of $G$. For $N\unlhd G$ and $\theta\in\Irr(N)$, we write $\Irr(G|\theta)$ for the set of irreducible characters of $G$ lying over $\theta$. Similarly, write $\dz_p(G|\theta)$ for the set of $p$-defect-zero characters of $G$ lying over $\theta$.

Let $(R,\varphi)$ be a $p$-weight of $G$. Then $R=\mathrm{O}_p(\mathrm{N}_G(R))$ is called a \textit{$p$-radical subgroup} of $G$, and $\varphi \in \dz_p(\mathrm{N}_G(R)/R)$ is called \textit{the weight character}. In the following, when mentioning a weight, we often mean this version.  Note that $\varphi^\circ$, the restriction of \(\varphi\) to the set of all \(p'\)-elements of \(\mathrm{N}_G(R)/R\), is an irreducible projective Brauer character of $\mathrm{N}_G(R)/R$.
In this way, we can identify $(R,\varphi)$  with $(R,\varphi^\circ)$. We can also view $\varphi^\circ$ as an element of $\IBr_p(\mathrm{N}_G(R))$ via the inflation $\Inf$. Write $\overline{(R,\varphi)}$ for the $G$-conjugacy class containing $(R,\varphi)$ and $\WW_p(G)$ for the set of all such $\overline{(R,\varphi)}$. 

If a group $A$ acts on a finite group $G$ by automorphisms, then there is a natural action of $A$ on $\Irr(G)\sqcup \IBr_p(G)$ given by $\chi^a(g)=\chi(g^{a^{-1}})$ for any $g\in G$, $a\in A$ and $\chi\in\Irr(G)\sqcup \IBr_p(G)$. 
When a group $A$ acts on a finite set $X$, we write \(A_x\) for the stabilizer of \(x\in X\) in $A$, and  $x$ is said to be 
\textit{$A$-invariant} if $A_x=A$. If a group $A$ acts (from the right) on the sets $X$ and $Y$, then a map $f:X\longrightarrow Y$ is called \textit{$A$-equivalent} when $f(x)^a= f(x^a)$ for any $x\in X$ and $a\in A$.

\subsection{The Navarro-Alperin weight conjecture and a criterion}
Fix a prime $p$. Choose a positive integer $n$ sufficiently divisible and let $\mathcal K=\QQ_p(\zeta_n)$, its valuation ring $\mathcal O$ and its residue field $F$.

The triple $(\mathcal K,\mathcal O,F)$ is a splitting $p$-modular system for all finite groups under consideration, as in
\mcite{FFZ}{Section 1} or \cite[Section 2]{T2013}.
In the following, we write  $\HH_p=\operatorname{Gal}(\mathcal K/\QQ_p)$ and $\zeta_n$ for a chosen primitive $n$-th root of
unity in a fixed algebraic closure of $\QQ_p$.
Let $K$ be the subfield of $\mathcal K$, generated over $\QQ$
by all roots of unity in $\mathcal K$ whose orders are prime to $p$. The image of the restriction map $\HH_p \rightarrow \Gal(K/\QQ)$ is the cyclic group generated by the automorphism $\sigma_p$ defined by $\zeta^{\sigma_p}=\zeta^p$ for every root of unity $\zeta\in K$ of order prime to $p$.

All ordinary characters are regarded as $\mathcal K$-valued.
Brauer characters are defined on $p$-regular elements and are
also regarded as $\mathcal K$-valued, using the lifts to
$\mathcal O$ of roots of unity in $F$ of order prime to $p$.
The action of $\HH_p$ on ordinary and Brauer characters is given by first evaluating the character and then applying the Galois automorphisms to the character values.
The action of $\HH_p$ on ordinary characters and Brauer characters is as follows by first compute the character, then apply the Galois automorphism to values of the character.

Write $\Aut(G)$ for the automorphism group of $G$.  For any \(a=(\sigma,\alpha)\in\HH_p\times\Aut(G)\), and for any subgroup \(H\le G\) and \(\chi\in \IBr_p(H)\sqcup\Irr(H)\), we define \(H^a=H^\alpha\) and \(\chi^a\in \IBr_p(H^a)\sqcup\Irr(H^a)\) by 
\(\chi^a(h)=\chi(h^{\alpha^{-1}})^\sigma\) for any 
$h\in H^a$. The group \(\HH_p\times\Aut(G)\) acts on the set of $p$-weights \((R,\varphi)\) of $G$ by the automorphism 
\(\alpha\in\Aut(G)\) sending  $R$ to $R^\alpha$ and transports the weight character \(\varphi\) along the induced isomorphism 
\(\NN_G(R)/R\rightarrow\NN_G(R^\alpha)/R^\alpha\), while \(\sigma\in\HH_p\) acts on the resulting character. This action commutes with $G$-conjugation, hence induces an action on the set \(\WW_p(G)\).

With these notations as above, we state the Navarro-Alperin Weight conjecture from \cite{FFZ}.
\begin{conjecture}[{NAW conjecture}]\label{con2.1}
Let $G$ be a finite group and $\HH_p$  described as above. Then there exists an $(\HH_p\times \Aut(G))$-equivalent bijection from $\IBr_p(G)$ to $\WW_p(G)$.    
\end{conjecture}

We want to show that a finite non-abelian simple group  satisfies the inductive NAW condition at the prime $p$ from \cite[Definition 7.1]{FFZ}. Feng, Fu and Zhou in \cite[Proposition 9.1]{FFZ} have given a technical criterion of this inductive condition. Now we only state a its simplified version that implies this inductive  condition for simple groups with trivial Schur multiplier and cyclic outer automorphism group.
\begin{proposition}\label{prop:criterion}
Let $G$ be a finite non-abelian simple group with trivial Schur multiplier and cyclic outer automorphism group. 
Assume that there exists an $(\mathcal H_p\times\Aut(G))$-equivariant bijection
\begin{equation}\label{equ:bijection}
 \Omega:\IBr_p(G)\longrightarrow\mathcal W_p(G).
\end{equation}
For $\eta\in\IBr_p(G)$, write
$\Omega(\eta)=\overline{(R,\varphi)}$.
Assume that there are finite groups $A$, $\widetilde A$ with $G\unlhd A\unlhd\widetilde A$ such that, under the conjugation homomorphism $\widetilde A\to\Aut(G)$ with kernel $\CC_{\widetilde A}(G)$, we have
\[
\widetilde A/\CC_{\widetilde A}(G)\cong\Aut(G)_{\eta^{\HH_p}},\quad A/\CC_{\widetilde A}(G)\cong\Aut(G)_\eta,
\]
and $\CC_{\widetilde A}(G)=\ZZ(A)$.
Here, $\Aut(G)_{\eta^{\HH_p}}:=\{\alpha\in\Aut(G)\mid \eta^\alpha\in \eta^{\HH_p}\}$ denotes the stabilizer of the $\HH_p$-orbit of $\eta$ in $\Aut(G)$.
Assume that  $\eta$ and $\varphi^\circ$ have extensions $\widehat\eta\in\IBr_p(A)$ and 
$\widehat\varphi\in\IBr_p(\NN_A(R))$ lying over the same central character of $\ZZ(A)$, respectively, and that
\begin{equation}\label{extension criterion}
\text{$\widehat \eta$ and $\widehat \varphi$ are invariant under the stabilizer $(\HH_p\times\NN_{\widetilde A}(R))_{\varphi}$}.
\end{equation}
Then $G$ satisfies the inductive NAW condition at the prime $p$.
\end{proposition}

\begin{remark}\label{rem: 2.2}
The above Proposition can be simplified. For example, suppose that $\widetilde{A}\le\Aut(G)$, where $G$ is identified with the inner automorphism group $\Inn(G)$. If $\CC_{\widetilde A}(G)=\ZZ(A)=1$, and $A/G$ is a $p'$-group, then the condition \eqref{extension criterion} is automatically satisfied, as explained in \cite[Remark 9.2]{FFZ}.
\end{remark}

The following results are often used to verify  the the condition \eqref{extension criterion}.

\begin{lemma}\label{lem:h3}
Assume that $G \unlhd A$ and that $A/G$ is cyclic.
Let $\theta\in\IBr_p(G)\sqcup \Irr(G)$ be $\HH_p$-invariant. If  $(|A_\theta/G|,p-1)=1$, then $\theta$ has a unique $\HH_p$-invariant extension to $A_\theta$.
\end{lemma}
\begin{proof}
Cyclic quotient group shows that  $\theta$ extends to $A_\theta$. 
Choose its an extension $\chi\in\IBr_p(A_\theta)$ and let $X=\IBr_p(A_\theta/G)$ the group of Brauer characters of $A_\theta/G$. By \cite[Corollary 8.20]{N1998}, any extension of $\theta$ to $A_\theta$ is of the form $\lambda\chi$ for some $\lambda\in X$. 

For $\sigma_p\in\HH_p$, since $\theta$ is $\sigma_p$-invariant, there exists $\mu\in X$ such that $\chi^{\sigma_p}=\mu\chi$. For $\lambda\in X$, we have $(\lambda\chi)^{\sigma_p}=\lambda^p \mu\chi$. Thus $\lambda\chi$ is $\sigma_p$-invariant if and only if $\lambda^p\mu = \lambda$, i.e., $\lambda^{p-1}=\mu^{-1}$.
The group $X$ is finite abelian, and $(|X|,p-1)=1$ by hypothesis.
Thus the $(p-1)$-th power map on $X$ is  bijective, and there is a unique $\lambda_0\in X$ satisfying $\lambda_0^{p-1}=\mu^{-1}$.
Therefore, $\chi_0=\lambda_0\chi$ is the unique $\sigma_p$-invariant extension of $\theta$ to $A_\theta$. 
For ordinary characters, using $X=\Irr(A_\theta/G)$ and \cite[Corollary 6.17]{I1976}, and we are done.
\end{proof}
\begin{remark}
The hypothesis $(|A_\theta/G|,p-1)=1$ cannot be dropped.
For instance, take $p=7$, $G=C_3\lhd A=C_9$, so $A/G\cong C_3$.
Let $\theta\in\Irr(G)$ be non-trivial. Since $7\equiv 1\pmod{3}$, $\theta$ is $\sigma_7$-fixed. The three extensions of $\theta$ to $A$ are permuted in a 3-cycle by $\sigma_7$, so there is no $\sigma_7$-invariant extension. In fact, we use Lemma \ref{lem:h3} in this paper only at $p=3$.
\end{remark}

\begin{lemma}\label{lem:real}
Assume that $G\unlhd A$ and that $A/G$ has odd order.
If $\theta\in\IBr_p(G)\sqcup\Irr(G)$ is real-valued, then it has a unique  $(A\times\HH_p)_\theta$-invariant real-valued extension to $A_\theta$. 
\end{lemma}
\begin{proof}
For ordinary characters, this is \cite[Lemma 4.2]{Johansson}.
For $\theta\in\IBr_p(G)$, similar to the proof of Lemma \ref{lem:h3}, we only need to replace $\chi^{\sigma_p}$ with $\overline{\chi}$. We can still obtain  that $\theta$ has a unique real-valued extension.
Moreover, the uniqueness ensures that the extension is also 
$(A\times\HH_p)_\theta$-invariant.
\end{proof}

\subsection{Some notation and conventions for the Suzuki and Ree groups}\label{section 2.3}
Let \(G\) be a Suzuki group \({}^{2}B_{2}(2^{2m+1})\) for
\(m\geq 1\), or one of the Ree groups
\({}^{2}G_{2}(3^{2m+1})\) and \({}^{2}F_{4}(2^{2m+1})\)
for \(m\geq 0\). Let \(\ell\) denote its defining characteristic, namely, \(\ell=2\) for \({}^{2}B_{2}\) and \({}^{2}F_{4}\), and \(\ell=3\) for \({}^{2}G_{2}\).
We write $q^{2}:=\ell^{2m+1}$. Let \(p\) be a prime divisor of \(|G|\) with \(p\neq\ell\). Except for \({}^{2}B_{2}(8)\), \({}^{2}G_{2}(3)\), and \({}^{2}F_{4}(2)\), the groups in these families are non-abelian simple groups with trivial Schur multiplier. More precisely, \({}^{2}B_{2}(8)\) is simple but has Schur multiplier \(C_{2}\times C_{2}\), whereas
\({}^{2}G_{2}(3)\) and \({}^{2}F_{4}(2)\) are not simple.

Let \(\mathbf G\) be a connected simple algebraic group of
adjoint type \(B_{2}\), \(G_{2}\), or \(F_{4}\), respectively,
defined over an algebraic closure of \(\mathbb F_{\ell}\).
Write \(F_{\ell}\) for the standard Frobenius endomorphism of
\(\mathbf G\). Let \(\gamma\) be the exceptional graph
isogeny associated with the corresponding Dynkin type, satisfying $\gamma^{2}=F_{\ell}$ and $\gamma\circ F_{\ell}=F_{\ell}\circ\gamma$. Define the Steinberg endomorphism
$F=F_{\ell}^{m}\circ\gamma$. Then the fixed-point group $\mathbf G^{F}\cong G$. The outer automorphism group \(\Out(G)\) of \(G\) is cyclic of order \(2m+1\), generated by the field automorphism induced by \(F_\ell\). Then 
$\Aut(G)\cong G\rtimes \langle F_{\ell}\rangle$.
Furthermore, \((\mathbf G,F)\) is self-dual and 
the center of $\bf G$ is trivial (see \cite[Theorem 2.5.12]{G1998}).

\section{The simple groups with cyclic Sylow subgroups}\label{section 3}
In the Section, we keep the notations in Subsection \ref{section 2.3} and consider the simple Suzuki and Ree groups with cyclic Sylow subgroup. 
\begin{lemma}\label{cyclic Sylow}
Let $G$ be a finite non-abelian simple group of order divisible by $p$ and $L$ the universal $p'$-covering group of $G$. Suppose that there exists a finite group $A$ such that $G\unlhd A$ and $A$ induces all automorphisms on $G$ by conjugation. 
If $A/G$ is cyclic, $\CC_A(G)=\ZZ(A)$ and $G$ has a cyclic Sylow $p$-subgroup, then $G$ satisfies the inductive NAW condition at all  primes $p$.
\end{lemma}
\begin{proof}
See \cite[Corollary 10.4]{FFZ}.
\end{proof}

In the following, we use the notion of Sylow \(\Phi^{(p)}\)-tori of $\bf G$ as in \cite[Section 8.1]{G2007}. Therefore, let \(\Phi^{(p)}\) be the cyclotomic polynomial over \(\mathbb{Q}(\sqrt{\ell})\) that divides the generic order of \(G\) and satisfies \(p \mid \Phi^{(p)}(q)\). This is well-defined except for \(p= 2\) and \(G={}^2G_2(q^2)\), or \(p= 3\) and \(G={}^2F_4(q^2)\). Using this definition, we obtain the following facts.

\begin{theorem}\label{thm:Suzuki}
Let $G={}^2B_2(q^2)$ and $p$ be odd. Then $G$ satisfies the inductive NAW condition at the prime $p$.
\end{theorem}
\begin{proof}
Let $G={}^2B_2(q^2)$ with $q^2=2^{2m+1}$ and $m\geq 1$. Then $|G|=q^4(q^2-1)(q^4+1)$. Note that $p$ divides $q^2-1$ or $q^4+1=(q^2+\sqrt{2}q+1)(q^2-\sqrt{2}q+1)$. By \cite[Section 16]{I2007}, any Sylow $p$-subgroup of $G$ is contained in a cyclic Sylow $\Phi^{(p)}$-torus of $\bf G$.  Applying Lemma \ref{cyclic Sylow} to $A=\Aut(G)$, and we complete the proof.
\end{proof}

\begin{proposition}\label{lem:3.2}
Let $G={ }^2F_2(q^2)$ with $q^2=2^{2m+1}$ and $m\geq 1$. If $p$ is odd and  $p\nmid q^8-1$, then $G$ satisfies the inductive NAW condition at the prime $p$.
\end{proposition}
\begin{proof}
Note that the order of $|G|$ is $q^{24}(q^2-1)^2(q^2+1)^2(q^4+1)^2(q^4-q^2+1)(q^8-q^4+1)$.
Since $p$ is odd and $p\nmid q^8-1=(q^4+1)(q^2+1)(q^2-1)$, 
we have \(p\neq 3\), and thus \(p\ge 5\). Moreover, $p$ divides either \(q^4-q^2+1\) or \(q^8-q^4+1\). In this case, any Sylow $p$-subgroup of $G$ is contained in a Sylow $\Phi^{(p)}$-torus of $\bf G$. By \cite[TABLE I]{An}, there exists such cyclic Sylow $\Phi^{(p)}$-torus of order $q^4-q^2+1$. Also, $q^8-q^4+1=u_+u_-$, where $u_\pm=q^4+q^2+1\pm \sqrt{2}q(q^2+1)$. Similarly, by \cite[TABLE I]{An}, both $u_\pm$ occur as orders of cyclic Sylow tori. Their greatest common divisor in $\QQ(\sqrt{2})$ is $1$ because a common odd divisor would divide both $q^4+q^2+1$ and $q^2+1$.
In particular, any odd prime $p$ divides at most one of 
\(u_+,u_-\). Thus there exists a cyclic torus contains a Sylow $p$-subgroup in every case under consideration. By Lemma \ref{cyclic Sylow} to $A=\Aut(G)$, and we are done. 
\end{proof}

\section{The Ree groups \texorpdfstring{${}^2F_2(2^{2m+1})$}{2F4}}\label{section 4}
In the Section, we will verify all hypothesis of Proposition \ref{prop:criterion} for the simple Ree groups $G={ }^2F_2(q^2)$ with $q^2=2^{2m+1}$ and $m\geq 1$. The Tits group ${}^2F_4(2)$ is handled in Section \ref{section 4}.  
By Proposition \ref{lem:3.2}, it suffices to consider the following two cases in Subsections \ref{p>=5} and \ref{p=3}.
\begin{equation}\label{eq:remaining}
 p=3,\ \text{or}\ p\geq5,\ p\mid(q^2-1)(q^2+1)(q^2-\sqrt{2}q+1)(q^2+\sqrt{2}q+1).
\end{equation}

For $p\geq 5$, by \cite[Corollary~4.3 and  Appendix A]{Himstedt} and together with \cite{GH}, the ordinary basic set $\mathcal{B}$, which is the union of these 
Lusztig series $\mathcal{E}(G,s)$,  is unitriangular,
where $s$ runs over the set of semisimple $p'$-elements of $G$.
Write $D:=\Out(G)\cong\langle F_2\rangle$. It is clear that the set $\mathcal{B}$ is invariant under the action of $\HH_p\times D$. Therefore, by \cite[Lemma~10.7]{FFZ}, there exists  an $(\HH_p\times D)$-equivalent bijection 
\begin{equation}\label{eq:basic}
\mathcal{B}\longrightarrow\IBr_p(G).
\end{equation}
We will construct an $(\HH_p\times D)$-equivalent bijection between $\mathcal{B}$ and $\WW_p(G)$ in Subsection \ref{p>=5}.

For $p=3$, which is a bad prime for $\bf G$, the results on ordinary basic sets in \cite{GH} do not apply in this
situation. For example, the set $\mathcal{E}(G,1)$ of ordinary unipotent irreducible characters of $G$ is not a basic set(see \cite[Section 4]{Himstedt}). However, by \cite[Appendix C]{Himstedt}, we can still select a set of characters to form a unitriangular basic set(see details in \eqref{principal 3-block} below). The situation for the non-unipotent characters of $G$ is less complicated. Using the character table of $G$ in \cite{Himstedt}, one can see that for all non-trivial semisimple $3'$-elements $s$, the set $\bigcup_{1\neq s}\mathcal{E}(G,s)$ is an ordinary basic set for a union of blocks. The corresponding decomposition numbers are given in \cite[Appendix D]{Himstedt}.

\subsection{Equivalent bijection for \texorpdfstring{$p\geq 5$}{p>=5}}\label{p>=5}
In the subsection, we assume that  $p\geq5$ and $p\mid(q^2-1)(q^2+1)(q^2-\sqrt{2}q+1)(q^2+\sqrt{2}q+1)$.
Then it follows from \cite{Johansson} that any Sylow $p$-subgroup  of $G$ is abelian of rank-two. Based on the different of $p$-radical subgroups, we will discuss them separately.

\begin{lemma}[{\cite[(2B)]{An}}]\label{p3:prop:radical}
There are exactly two non-trivial \(G\)-conjugacy classes of \(p\)-radical subgroups, represented by $R$ and $P$, where $R$ is cyclic and $P$ is non-cyclic abelian.
\end{lemma}
\subsubsection{Non-cyclic abelian \texorpdfstring{$p$}{p}-radical subgroups}
Let $P$ be non-cyclic abelian \(p\)-radical subgroup of $G$. Then by \cite[(2B)-(b)]{An}, $P$ is a Sylow $p$-subgroup of $G$. Moreover, $\CC_G(P)=T$ and $P=\OO_p(T)$, where $T={\bf T}^F$ and ${\bf T}$ is a Sylow $\Phi^{(p)}$-torus of $\bf G$. 
The normalizer of the torus ${\bf T}$ as following Table \ref{tab:1} (\cite[TABLE I]{An}). It is not difficult to verify that $\NN_G(P)=\NN_G({\bf T})\cong T\rtimes W$, where $W$ is a $p'$-group.
\begin{table}[ht]
\centering
\caption{Sylow tori and normalizers}\label{tab:1}
\begin{tabular}{@{}lllrr@{}}
\toprule
Prime condition &Tori $T$ & $W=\NN_G({\bf T})/T$\\
\midrule
$p\mid q^2-1$ & $C_{q^2-1}\times C_{q^2-1}$ & $D_{16}$\\
$p\mid q^2+1$ & $C_{q^2+1}\times C_{q^2+1}$ & $\GL_2(3)$\\
$p\mid q^2\pm\sqrt{2}q+1$ & $C_{q^2\pm\sqrt{2}q+1}\times C_{q^2\pm\sqrt{2}q+1}$ & $G_8$\\
\bottomrule
\end{tabular}
\end{table}
Here, $D_{16}$ is a dihedral group of order $16$, and $G_8$ denotes the Shephard-Todd complex reflection group of order $96$. In this case, we have the following result. Write $\Irr_{p'}(G)$ for the set of all $p'$-characters(their degrees are not divisible by $p$) of $G$.
\begin{proposition}\label{bij: rank-two}
Let $P$ be a Sylow $p$-subgroup of $G$. Then there exists an $(\HH_p\times \Aut(G)_P)$-equivalent bijection
$\Omega_P:\Irr_{p'}(G)\cap \mathcal{B}\longrightarrow\dz_p(\NN_G(P)/P)$.
\end{proposition}
\begin{proof}
We identify $G$ and its subgroups with their dual groups. 
Fix the standard isomorphism $T\cong\Irr(T), s\mapsto\theta_s$. For any \(s\in T\), let $N=\NN_G(P)$, $N_s=\NN_G(P)_{\theta_s}$ and $W_s=\NN_G({\bf T})_{\theta_s}/T$. Through the dual identification, \(W_s\) is identified with the relative Weyl group associated with the chosen torus inside the centralizer \(\CC_{G}(s)\).

Note that $N=\NN_G({\bf T})$, by \cite[Lemma 3.3]{Johansson}, we have an $(\HH_p\times \Aut(G)_P)$-equivalent  extension map 
\[
\Lambda:\Irr(T)\longrightarrow\bigcup_{T\leq I\leq N}\Irr(I)
\]
with $\Lambda(\theta_s)\in\Irr(N_s)$ and $\Lambda(\theta_s)|_T=\theta_s$.
For any $\lambda\in\Irr(W_s)$, we define 
\[
\psi_{s,\lambda}=\Ind_{N_s}^{N}(\Lambda(\theta_s)\Inf_{W_s}^{N_s}(\lambda)).
\]
Since \(p\nmid |W|\), the degrees $\psi_{s,\lambda}(1)=|N:N_s|\lambda(1)$ are prime to \(p\). Consequently, we have $\psi_{s,\lambda}\in\Irr_{p'}(N)$. As in \cite[Section 3]{Johansson}, this gives all $p'$-irreducible characters of \(N\). For any $\sigma\alpha\in\HH_p\times \Aut(G)_P$, where $\sigma$ maps $p'$-th root of unity $\zeta$ to $\zeta^{p^j}$. By \cite[Proposition 3.4]{Johansson} and its proof, we obtain   $(\psi_{s,\lambda})^{\sigma\alpha}=
\psi_{\alpha(s^{p^j}),\ (\lambda^{\sigma})^{\alpha}}$.

We claim that $P\leq\ker(\psi_{s,\lambda})$ if and only if $s$ is $p'$-element. Note that \(\Inf_{W_s}^{N_s}(\lambda)\) is  equal to \(\lambda(1)\) on \(T\) and \(\Lambda(\theta_s)|_T=\theta_s\). Let $\Psi=\Lambda(\theta_s)\Inf_{W_s}^{N_s}(\lambda)$.
Then \(\Psi|_T=\lambda(1)\,\theta_s\).
By Mackey restriction(see \cite[Problem (5.6)]{I1976}),
we obtain
\[
\psi_{s,\lambda}|_T
=(\Ind_{N_s}^{N}(\Psi))|_T
=\sum_{x\in N\backslash N/N_s} \Ind_{T\cap (N_s)^{x}}^{T}(\Psi^x|_{T\cap (N_s)^x}).
\]
Since \(T\unlhd N\), we have \(T\cap (N_s)^x=T\), and then 
\((\Psi^x)|_T=(\lambda(1)\theta_s)^x=\lambda(1)(\theta_s)^x\).
Consequently, $\psi_{s,\lambda}|_T=\lambda(1)
\sum_{\theta_i}\theta_i$, where $\theta_i$ runs over the set of \(N\)-conjugate of \(\theta_s\). Since \(P\unlhd N\), all $\theta_i$ are trivial on \(P\) if \(\theta_s\) is trivial on \(P\). Conversely, if \(\psi_{s,\lambda}|_P\) is a multiple of \(1_P\), then every $\theta_i$ must be
trivial on \(P\). Therefore, $P\leq\ker(\psi_{s,\lambda})$ if and only if  $P\leq\ker(\theta_s)$. Since $P=\OO_p(T)$,
any  \(\theta_s\in\Irr(T)\) is trivial on \(P\) if and only if its order is prime to \(p\), that is, the corresponding $s$ is $p'$-element, as claimed.

Let $J_s:\mathcal{E}(G,s)\longrightarrow\mathcal{E}(\CC_{G}(s),1)$ be the Jordan decomposition of characters and 
$\mathcal I_s:\Irr(W_s)\longrightarrow
\mathcal{E}(\CC_{G}(s),{\bf T},1)$ the relative Weyl group correspondence(see \cite[Section 3]{Johansson}).  Define $\chi_{s,\lambda}=J_s^{-1}(\mathcal I_s(\lambda))$. Using \cite[Theorem B and Proposition~3.2]{Johansson}, for the above $\sigma\alpha\in\HH_p\times \Aut(G)_P$,  we have  also $(\chi_{s,\lambda})^{\sigma\alpha}=
\chi_{\alpha(s^{p^j}),\,(\lambda^\alpha)^{\sigma}}$.

Consider all the parametrization classes above  \(\Irr_{p'}(G)\). Restricting to semisimple  \(p'\)-elements  gives exactly the subset $\Irr_{p'}(G)\cap \mathcal{B}$ on global side. The same parametrization give exactly \(\Irr_{p'}(N/P)\) on the local side.  Since \(P\) is a Sylow \(p\)-subgroup of \(N\),  \(N/P\) is a \(p'\)-group. Thus $\dz_p(N/P)=\Irr_{p'}(N/P)$. Let \(\bar\psi_{s,\lambda}\) denote the character of \(N/P\) whose inflation is \(\psi_{s,\lambda}\), and define $\Omega_P(\chi_{s,\lambda})=\overline\psi_{s,\lambda}$. Then by the above process, $\Omega_P$ is an $(\HH_p\times \Aut(G)_P)$-equivalent bijection.
\end{proof}

\subsubsection{Cyclic \texorpdfstring{$p$}{p}-radical subgroups} 
Let $\mathbf{M}$ be a $\Phi$-split Levi subgroup of $\bf G$. Write $M={\bf M}^F$ and $Z=\ZZ(M)$. Let $R$ be a non-trivial cyclic $p$-radical subgroup of $G$. Then as in \cite[Section 3]{Malle}, $R=Z_p$, $M=\CC_{G}(R)$ and $\NN_{G}(R)=\NN_{G}({\bf M})$. 
In fact,  $R=\OO_p(\CC_G(R))$ by \cite[(2B)-(a)]{An},
and $M=Z\times K$. Moreover, by \cite[TABLE I]{An}, if $p|q^2\pm 1$, then $\NN_G(R)/\CC_G(R)\cong C_2$. If  $p|q^2\pm\sqrt{2}q+1$, then $\NN_G(R)/\CC_G(R)\cong C_4$.
The centralizer and normalizer of $R$ as in following Table 2. 
\begin{table}[ht]
\centering
\caption{Centralizers and normalizers of cyclic $p$-radical subgroups}\label{tab:2}
\begin{tabular}{@{}lllrr@{}}
\toprule
Prime condition & $\CC_G(R)=Z\times K$ & $\NN_G(R)$\\
\midrule
$p\mid q^2-1$ & $C_{q^2-1}\times\SL_2(q^2)$ & $(C_{q^2-1}\rtimes C_2)\times\SL_2(q^2)$\\
$p\mid q^2-1$ & $C_{q^2-1}\times{}^2B_2(q^2)$ & $(C_{q^2-1}\rtimes C_2)\times{}^2B_2(q^2)$\\
$p\mid q^2+1$ & $C_{q^2+1}\times\SL_2(q^2)$ & $(C_{q^2+1}\rtimes C_2)\times\SL_2(q^2)$\\
$p\mid q^2-\sqrt{2}q+1$ & $C_{q^2-\sqrt{2}q+1}\times{}^2B_2(q^2)$ & $(C_{q^2-\sqrt{2}q+1}\rtimes C_4)\times{}^2B_2(q^2)$\\
$p\mid q^2+\sqrt{2}q+1$ & $C_{q^2+\sqrt{2}q+1}\times{}^2B_2(q^2)$ & $(C_{q^2+\sqrt{2}q+1}\rtimes C_4)\times{}^2B_2(q^2)$\\
\bottomrule
\end{tabular}
\end{table}
Here, these structures are determined by the classification of centralizers of semisimple elements and the structure of  normalizers of tori (see \cite[TABLE I]{An} or \cite[Tables~B.9--B.10]{HH}). For example, we use the notation in \cite[TABLE I]{An}, from the first row to the end, corresponding to the centralizer $\CC_{G}(t_i)$($i=2,1,5,7,9$) contains the factor $A_1$, ${}^2B_2$, $A_1$, ${}^2B_2$, ${}^2B_2$, respectively. 

For completeness, we need to show that there is no proper radical subgroup whose centralizer is a maximal torus. To see this, suppose that $\CC_G(R)=T$ is a maximal torus. Then $T_p$ is normal in $\NN_G(R)$, so the radicality forces $R=T_p$, which is a Sylow subgroup. This is contradictory. 
Note that the exceptional semisimple centralizer $\CC_G(t_4)$ of type ${}^2A_2$ in \cite[TABLE I]{An} occurs at an element of order $3$, and does not add a $p$-radical subgroup for $p\geq5$.

Write $E=\NN_G(R)/\CC_G(R)$. Then $E\cong C_2$ or $C_4$. 
Since \(Z\) is cyclic, $Z=R\times Z_{p'}$,
and then
\[
\NN_G(R)/R=(Z_{p'}\rtimes E)\times K.
\]
As in \cite[Sections 3 and 5]{Malle}, by Clifford theory, every element of $\dz_p(\NN_G(R)/R)$ is parametrized by a triple \((\theta,\kappa,\epsilon)\), where \(\theta\in\Irr(Z_{p'})\), \(\kappa\in\dz_p(K)\) and \(\epsilon\in\Irr(E_\theta)\). In deed, the character \(\theta\) can be extended to $E_\theta$ by defining $\widehat\theta(ze)=\theta(z)$ for any $z\in Z_{p'}$ and $e\in E_\theta$. Let
\begin{equation}\label{eq:rankone}
 \psi_{\theta,\kappa,\epsilon}=
 \Ind_{(Z_{p'}\rtimes E_\theta)\times K}^{(Z_{p'}\rtimes E)\times K}
((\widehat\theta\epsilon)\otimes\kappa).
\end{equation}
Note that $\psi_{\theta,\kappa,\epsilon}(1)
=|E:E_\theta|\kappa(1)$. Since \(|E:E_\theta|\) is a \(p'\)-number, we have $\psi_{\theta,\kappa,\epsilon}(1)_p=\kappa(1)_p$.
On the other hand, $|\NN_G(R)/R|_p=|K|_p$.
Thus $\psi_{\theta,\kappa,\epsilon}\in\dz_p(\NN_G(R)/R)$ if and only if  $\kappa\in\dz_p(K)$.
\begin{lemma}\label{E}
$E_\theta=E$ for $\theta=1$, and $E_\theta=1$ otherwise.    
\end{lemma}
\begin{proof}
Suppose that  \(E=C_2\). Then the action of \(C_2\) on \(Z_{p'}\) is inversion. Since \(|Z_{p'}|\) is odd, $E_\theta=C_2$ for $\theta=1$ and $E_\theta=1$ otherwise. Suppose that  \(E=C_4\), the square of the action  is inversion. Accordingly, we may choose a generator so that the action is given by the \(q^2\)-th power.
Note that $(q^2+\sqrt{2}q+1)(q^2-\sqrt{2}q+1)=q^4+1$  and $q^4\equiv-1\pmod{q^2\pm\sqrt{2}q+1}$.
If a non-trivial subgroup stabilizes \(\theta\), then it contains the unique involution of \(C_4\), and so we again have \(\theta=\theta^{-1}\). Thus
$E_\theta=C_4$ for $\theta=1$ and $E_\theta=1$ otherwise.
\end{proof}

We can choose an odd-order cyclic subgroup $D_0\leq D$ representing the field automorphism such that it normalizes $Z_{p'}$, $K$ and satisfies $[D_0,E]=1$ as follows. Write $I=Z_{p'} \rtimes E$ and take a lift $D_1$ of the field-action preserving $I$. Since $E$ is the Sylow $2$-subgroup of $I$, the Frattini argument yields $I \rtimes D_1=I\NN_{I\rtimes D_1}(E)$.
Note that $\NN_I(E)=\CC_{Z_{p'}}(E)\rtimes E$.
Since $E=1,\ C_2$ or $C_4$, and $Z_{p'}$ has odd order, $\CC_{Z_{p'}}(E)=1$, and then $\NN_I(E)=E$. Thus
$$1 \longrightarrow E \longrightarrow \NN_{I\rtimes D_1}(E) \longrightarrow D_1 \longrightarrow 1.$$
By the Schur-Zassenhaus theorem, we may choose an odd-order complement $D_0$ that is isomorphic to $D_1$ inside this normalizer. Note that $\Aut(E)=1$ or $C_2$. So the conjugation action of $D_0$ on $E$ must be trivial. Therefore, $D_0$ acts 
on all sets of their characters.
\begin{lemma}\label{lem:3.5}
Let $\sigma\alpha\in\HH_p\times D_0$. Then $\psi_{\theta,\kappa,\epsilon}^{\sigma\alpha}=\psi_{\theta^{\sigma\alpha},\kappa^{\sigma\alpha},\epsilon^{\sigma}}$.
\end{lemma}
\begin{proof}
Since $|Z_{p'}|$ is odd and $E_\theta$ has order $1$, $2$ or $4$, we know that $\widehat \theta$ is the unique linear extension of odd degree. Since $E_\theta=E$ or $E_\theta=1$ by Lemma \ref{E}, $D_0$ acts trivially on $E_\theta$, and then any $\epsilon$ is $D_0$-fixed. For any $\sigma\in\HH_p$ and $\alpha\in D_0$, $\widehat \theta^{\sigma\alpha}$ is an extension of $\theta^{\sigma\alpha}$. By uniqueness of the extension, we have $\widehat \theta^{\sigma\alpha}=\widehat{ \theta^{\sigma\alpha}}$. Both induction and tensor product of characters commute with the action of $\sigma\alpha$. We complete the proof. 
\end{proof}

\begin{lemma}\label{defect zero characters of $K$}
Let $K$ be as in Table \ref{tab:2}. When $K=\SL_2(q^2)$, its unipotent characters are  $1, \St$, and when $K={}^2B_2(q^2)$, its unipotent characters are $1, \St, W^+, W^-$. 
Moreover, for any \(\sigma\alpha\in\HH_p\times D_0\),  
$\St^{\sigma\alpha}=\St$, $(W^\pm)^\alpha=W^\pm$ and 
\[(W^\pm)^\sigma=
		\begin{cases}
			W^\pm,&\zeta_4^\sigma=\zeta_4,\\
			W^\mp,&\zeta_4^\sigma=-\zeta_4.
		\end{cases}
\]
Here, $\St$ denotes the Steinberg character, $W^+$ and $W^-$ have the same degree.
\end{lemma}

\begin{proof}
Suppose that \(K=\SL_2(q^2)\). Then its two unipotent characters are $1$ and $\St$ with degrees $1$ and $q^4$(see \cite[Section 4]{W13}).
Suppose that $K={}^2B_2(q^2)$. Then there are four unipotent characters $1,\St,W^+,W^-$ with degrees $1,q^4,(\sqrt{q^2/2})(q^2-1),(\sqrt{q^2/2})(q^2-1)$ (see \cite[Section 2]{Ghaffarzadeh}).  

Clearly, $\St$ is fixed because Galois automorphisms and field automorphisms preserve degrees of characters. For \(W^+\) and \(W^-\), they are individually fixed by field automorphisms
(see \cite[Lemma 2.2]{Ghaffarzadeh}). By \cite[proof of Theorem B]{Johansson}, their Galois actions are determined by the action on \(\zeta_4\). They are fixed if \(\zeta_4^\sigma=\zeta_4\), and interchanged if \(\zeta_4^\sigma=-\zeta_4\). 
\end{proof}

Let $B$ be a $p$-block of $G$ parametrized by the cuspidal pair $({\bf M},\lambda)$ of $\bf G$, where ${\bf M}$ is $\Phi$-split
Levi subgroup of ${\bf G}$. Recall that the corresponding relative Weyl group $W_G({\bf M},\lambda)=\NN_{G}({\bf M},\lambda)/M$(see \cite[Section 3]{Malle}). Let $R$ be \(p\)-radical subgroup of $G$. Then by \cite[Corollary 3.2 and Remark 2.5]{Malle}, we have $R=Z_p=\OO_p(M)$, and then $M/R=Z_{p'}\times K$. 
\begin{proposition}\label{bij: rank-one}
Let $R$ be cyclic \(p\)-radical subgroup of $G$. Then 
any bijection between $\mathcal{E}(G,s)\cap\Irr(B)$ and $\dz_p(\NN_G(R)/R\mid\lambda)$ is $(\HH_p\times D)_{{\bf M},\lambda}$-equivalent.
\end{proposition}
\begin{proof}
For the cuspidal pair $({\bf M},\lambda)$, note that $\NN_{G}({\bf M})=\NN_G(R)$ and $M=\CC_G(R)$. So $W_G({\bf M},\lambda)\cong E_\lambda$. Thus the parameterization chain in \cite[Propositions~3.5-3.6 and Corollary~3.7]{Malle} is
\begin{equation}\label{eq:chain}
\mathcal{E}(G,s)\cap\Irr(B)\longleftrightarrow\Irr(E_\lambda)
 \longleftrightarrow\dz_p(\NN_G(R)/R\mid\lambda),\ 
\chi_{\lambda}\longleftrightarrow \epsilon\longleftrightarrow \psi_{\lambda,\epsilon}.
\end{equation}
Here, $\lambda=\theta\otimes\kappa\in\dz_p(M/R)$; $\psi_{\lambda,\epsilon}$ is identified with $\psi_{\theta,\kappa,\epsilon}$ in \eqref{eq:rankone}; The character $\chi_{\lambda}$ denotes an irreducible constituent of the generalized Harish-Chandra induction $R_{\bf M}^{\bf G}(\lambda)$. 

Note that $E_\lambda=1, C_2$ or $C_4$.  We next will show that $(\HH_p\times D)_{{\bf M},\lambda}$ acts trivially pointwise on both sets of \eqref{eq:chain}. If $E_\lambda=1$, then both sets in \eqref{eq:chain} are singletons. This is trivial. So we consider two cases of $E_\lambda=C_2$ and $E_\lambda=C_4$. Let $\sigma_{p}\alpha\in(\HH_p\times D)_{{\bf M},\lambda}$.
Then it is clear that $\alpha$ normalizes $Z_{p'}\times K$, and that $\alpha$ centralizes $E_\lambda$.

Suppose that $E_\lambda=C_2$. Then we  assume that $\Irr(E_\lambda)=\{1,sgn\}$, where $sgn$ is the symbol character. The corresponding local characters are $\psi_{\lambda, 1}$ and $\psi_{\lambda,sgn}$. Since  $1$ and $sgn$ are rational-valued, and together with Lemma \ref{lem:3.5}, $\psi_{\lambda, 1}$ and $\psi_{\lambda,sgn}$ are $\sigma_p\alpha$-fixed.
On the other hand, the corresponding two global characters, denoted by $\chi_{\lambda_\pm}\in\mathcal{E}(G,s)$. By Jordan decomposition of characters and Table \ref{tab:2}, they correspond to two unipotent characters of $\CC_G(s)$, where $\CC_{G}(s)$ is of type $A_1$ for $s\in\{t_2,t_5\}$ or ${}^2B_2$ for $s=t_1$(see Table \ref{tab:2} above). 
If $\CC_{G}(s)$ is of type $A_1$, then two $\chi_{\lambda_\pm}$
must correspond to $1, \St$, which  have  different degrees.
If $\CC_{G}(s)$ is of type ${}^2B_2$, then two $\chi_{\lambda_\pm}$ correspond to two unipotent characters with different degrees expect for $W^+$ and $W^-$. For this exceptional case, two $\chi_{\lambda_\pm}$ are unipotent, and they come from $R_{\bf M}^{\bf G}(W^\pm)$ induced from the two unipotent characters $W^\pm$ of $M=C_{q^2-1}\times {}^2B_2(q^2)$. Use the notation in \cite{Himstedt}, $\chi_{\lambda_-}$ and $\chi_{\lambda_+}$ come from 
the series $\{\chi_2,\chi_{19}\}$ or $\{\chi_3,\chi_{20}\}$. The two degrees within each series are indeed different by \cite[TABLE A.1]{Himstedt}. In either case, two $\chi_{\lambda_\pm}$ have different degrees. Therefore, they are $\sigma_p\alpha$-fixed. 

In the last case of the above paragraph, since $p|q^2-1$ and $p$ is odd, by \cite[Lemma 3.1(a)]{Johansson}, we have $p\equiv 1,7\pmod 8$.  Then $\zeta_4^{\sigma_p}=\zeta_4^{p}=\pm\zeta_4$. By Lemma \ref{defect zero characters of $K$}, $\sigma_p$ may interchange $W^+$ and $W^-$. When $(W^+)^{\sigma_p}=W^-$, $\sigma_p$  interchange  the two corresponding global series 
$\{\chi_2,\chi_{19}\}$, $\{\chi_3,\chi_{20}\}$ by $\chi_2\leftrightarrow \chi_3$ and $\chi_{19}\leftrightarrow \chi_{20}$, and $\psi_{W^+,\epsilon}^{\sigma_p}=\psi_{W^-,\epsilon}$ by Lemma \ref{lem:3.5}. However, this process does not change the stabilizer-actions on the global and local  characters.

Suppose that $E_\lambda=C_4$. Let $c$ generate $C_4$.
Then all characters of $C_4$ are $\eta_j(c)=\zeta_4^j$,
where $\zeta_4$ is a primitive fourth root of unity and 
$j=0,1,2,3$. Thus we write
$\Irr(E_\lambda)=\{\eta_j|j=0,1,2,3\}$. In this case, by Table \ref{tab:2}, \(p\mid q^4+1\) and \(p\) is odd. So by \cite[Lemma 3.1(d)]{Johansson}, $p\equiv 1\pmod 4$.  
Then $\zeta_4^{\sigma_p}=\zeta_4^{p}=\zeta_4$, and thus every $\eta_j$ is $\sigma_p$-fixed. Clearly, they are also $\alpha$-fixed. Therefore, by Lemma \ref{lem:3.5}, these characters $\psi_{\lambda,\eta_j}$ are $\sigma_p\alpha$-fixed. 
For the global side, the four characters, denoted by  $\chi_{\lambda_j}(j=0,1,2,3)$, correspond exactly to the four unipotent characters $1,\St,W^+,W^-$ of $\CC_G(s)$,
where $\CC_G(s)$ is of type ${}^2B_2(q^2)$ for $s\in\{t_7, t_9\}$(see  Table \ref{tab:2} above). Since $\zeta_4^{\sigma_p}=\zeta_4$, and together with Lemma \ref{defect zero characters of $K$}, all $1,\St,W^+,W^-$ are $\sigma_p\alpha$-fixed. Therefore, by Jordan decomposition of characters(see \cite[Theorem B]{Johansson}), these characters $\chi_{\lambda_j}(j=0,1,2,3)$ are $\sigma_p\alpha$-fixed.

Finally, we choose a representative $({\bf M},\lambda)$ from each $(\HH_p\times D)$-orbit. By the correspondence in \eqref{eq:chain}, both sets have the same number of elements. 
Moreover, from the above process, since $(\HH_p\times D)_{{\bf M},\lambda}$ acts trivially on both sets, any bijection from 
$\mathcal{E}(G,s)\cap\Irr(B)$ to $\dz_p(\NN_G(R)/R\mid\lambda)$
is automatically $(\HH_p\times D)_{{\bf M},\lambda}$-equivalent.
We complete the proof.
\end{proof}

\begin{theorem}\label{equ:bijection when p>=5}
Then there exists an $(\HH_p\times\Aut(G))$-equivalent bijection
$\Omega:\IBr_p(G)\longrightarrow\WW_p(G)$.
\end{theorem}
\begin{proof}
First, when the Sylow $p$-subgroup $P$ of $G$ as the first component of the weights, we apply Proposition \ref{bij: rank-two}, there exists an $(\HH_p\times \Aut(G)_P)$-equivalent bijection $\Omega_P:\Irr_{p'}(G)\cap\mathcal{B}\longrightarrow\dz_p(\NN_G(P)/P)$.
For $\chi\in \Irr_{p'}(G)\cap\mathcal{B}$, define $\widetilde\Omega_P(\chi)=\overline{(P,\Omega_P(\chi))}$.
For any \(\alpha\in\Aut(G)\), there is an inner automorphism \(c_g\) of $G$ such that $\beta=c_g\circ\alpha\in\Aut(G)_P$. Note that \(\alpha\) and \(\beta\) induce the same actions on these two sets. Thus $\widetilde \Omega_P$ is an $(\HH_p\times \Aut(G))$-equivalent bijection.

Second, when cyclic $p$-subgroups $R$ of $G$ as the first component of the weights, by Proposition \ref{bij: rank-one}, we choose $\Omega_R$ as an
$(\HH_p\times D)_{{\bf M},\lambda}$-equivalent bijection between
$\mathcal{E}(G,s)\cap\Irr(B)$ and $\dz_p(\NN_G(R)/R\mid\lambda)$. For $\sigma_p\alpha\in \HH_p\times D$ and $\chi\in \mathcal{E}(G,s)\cap\Irr(B)$, define $\widetilde\Omega_R(\chi^{\sigma_p\alpha})=\Omega_R(\chi)^{\sigma_p\alpha}$. It is easy to check that $\widetilde\Omega_R$ is
an $(\HH_p\times D)$-equivalent bijection. Note that $\Aut(G)=G\rtimes D$, so it is also $(\HH_p\times \Aut(G))$-equivalent.

Finally, for any $\chi\in\dz_p(G)$, the assignment $\chi\mapsto\overline{(1,\chi)}$ is obviously equivalent.
Apply the bijections $\widetilde\Omega_P$ and $\widetilde\Omega_R$, and together with the bijection \eqref{eq:basic}, we obtain an $(\HH_p\times\Aut(G))$-equivalent bijection between  $\IBr_p(G)$ and $\WW_p(G)$. We complete the proof.
\end{proof}

\subsection{Compatible extensions for \texorpdfstring{$p\geq5$}{p >= 5}}
Keep the groups $G$, $D$ and $\HH_p$ as above. We will verify that the bijection in Theorem \ref{equ:bijection when p>=5} satisfies the condition \eqref{extension criterion}. For global characters, since $\mathcal{B}$ is $(\HH_p\times \Aut(G))$-invariant, by \cite[Lemma 10.7(2)]{FFZ}, we only need to consider the ordinary characters in $\mathcal{B}$. Note that \(\Irr(G)\) contains a large number of real-valued characters, and only a portion of irreducible characters are not real-valued(see \cite{Himstedt}).

The real-valued characters are handled by Lemma \ref{lem:real}, and unipotent characters are handled by the following Lemma.

\begin{lemma}[{\cite[Proposition~4.3]{Johansson}}]\label{lem:unipotent}
If $\theta\in\Irr(G)$ is unipotent, then it has an $(\HH_p)_\theta$-invariant extension to $G\rtimes D$.  
\end{lemma}

For the remaining characters of $G$, with the CHEVIE numbering used in \cite[Table~A.2]{Himstedt} and \cite[Appendix~B]{HH}, the families of non-unipotent and non-real valued characters that need to be processed separately are
$\chi_{23},\chi_{24},\chi_{43},\chi_{44},\chi_{47},\chi_{48}$.
Moreover, the values of the last four contain non-zero imaginary parts $i$ on some conjugacy classes(see \cite[Table~B.12]{HH}).

\begin{lemma}\label{lem: The split Suzuki Levi families}
If $\chi$ is any of characters in the families 
$\chi_{23},\chi_{24}$, then it has an $(\HH_p\times D)_\chi$-invariant extension to $G\rtimes D_\chi$.  
\end{lemma}
\begin{proof}
Let ${\bf J=MU}$ be the standard parabolic subgroup of ${\bf G}$ with the split Levi subgroup ${\bf M}^F=M\cong C_{q^2-1}\times{}^2B_2(q^2)$.  For a non-trivial $\theta\in\Irr(C_{q^2-1})$ and a cuspidal unipotent character $W\in\{W^+,W^-\}$ as in Lemma \ref{defect zero characters of $K$}, we consider the Harish-Chandra induction
\begin{equation}\label{eq:hc}
 \chi=\Ind_{J}^{G}\Inf_{M}^{J}(\theta\otimes W).
\end{equation}
By \cite[Proposition~10.4(1) and its proof]{HHM}, $\chi$ is irreducible, and this gives the families $\chi_{23}$ or $\chi_{24}$.

First,  we claim that for any $d\in D$, there exists a $d$-fixed involution $n\in G$ such that 
\begin{equation}\label{eq: act}
(\theta\otimes W)^n=\theta^{-1}\otimes W.   
\end{equation}
To see this, we prove that there exists the element $n$ satisfies $z^n=z^{-1}$ for any $z\in C_{q^2-1}$ and $[n,{}^2B_2(q^2)]=1$. Let $w=w_0w_1$, where $w_0$ is the length element in root system of \({}^2F_4\) and the length element \(w_1\) in root subsystem $\Phi_1$ of type ${}^2B_2$ generated by two simple roots.  For a simple root $\alpha$, take the standard representative $n_\alpha = x_\alpha(1)x_{-\alpha}(-1)x_\alpha(1)$. These representatives satisfy $n_\alpha^2 = h_\alpha(-1)$ and the corresponding braid relations; their conjugate action by root groups has the form
$$\dot{w}x_\beta(t)\dot{w}^{-1}=x_{w\beta}(\varepsilon_{w,\beta}t),\quad \varepsilon_{w,\beta}\in\{1,-1\}.$$  
Since $-1=1$ in $\overline{\mathbb{F}}_2$, 
$n_\alpha^2=1$ and $ \varepsilon_{w,\beta}=1$.
Take $n=\dot{w}=\dot{w}_0\dot{w}_1$ and $n^2=1.$
Since $w\neq 1$, we have $n\neq 1$, and thus $n$ is an involution. For the above equation, we have 
$n x_\beta(t)n^{-1}=x_\beta(t)$ for $\beta\in\Phi_1$.  
These root groups generate the derived group of $M$, so $[n,{}^2B_2(q^2)]=1$. On the other hand, \(w_1\) acts trivially on $C_{q^2-1}$, while \(w_0=-1\) acts as inversion on the  maximal torus. Then \(w=w_0w_1\) also acts as inversion on $C_{q^2-1}$. Thus for any $z\in C_{q^2-1}$, we have 
$z^n=z^{-1}$.  Let $F=F_2^f\gamma$ the Steinberg endomorphism defining the Ree group, where $\gamma^2=F_2$ and \(\gamma\) is the exceptional automorphism. Then $\gamma(n)=n$, and so $F_2(n)=n$. Thus $n\in G$ and $n^d=n$ for any $d\in D$, as claimed.

Moreover, we claim that $J\cap J^n=M$. This fact will be used in Mackey's formula blew. Note that the root groups of $\bf J$ consist of  $\Phi_1$ and $\Phi^{+}\setminus\Phi_1^+$, where the latter generates ${\bf U}$. Since \(w_1\) preserves \(\Phi\setminus\Phi_1^+\), while \(w_0\) maps positive roots to negative roots, $w(\Phi^+\setminus\Phi_1^+) =\Phi^-\setminus\Phi_1^-$. Meanwhile, \(w\) fixes \(\Phi_1\) pointwise. Consequently, ${\bf J}^n=\bf J^-$ and $\bf J\cap\bf J^-=\bf M$, where \(\bf J^-\) is the opposite parabolic subgroup with the same Levi component $\bf M$. Taking \(F\)-fixed points yields $J\cap J^n=M$, as claimed.

Next, we will construct an extension of $\chi$ to $G\rtimes D_\chi$. Note that the character $W$ is $D_\chi$-fixed by Lemma \ref{defect zero characters of $K$}. So for any $d\in D_\chi$, Equalities  \eqref{eq:hc} and \eqref{eq: act} give $\theta^d=\theta^n\in\{\theta,\theta^{-1}\}$. Suppose that $\theta^d=\theta^{-1}$.
Let $d$ has odd order $r$ because $D_\chi$ has odd order. Then $\theta=\theta^{d^r}=\theta^{(-1)^r}=\theta^{-1}$, and this implies $\theta^2=1$. However, $C_{q^2-1}$ has odd order, which forces $\theta=1$, contradicting the assumption. So $\theta$ is $D_\chi$-fixed. Therefore, \(\theta\) can be extended to \(C_{q^2-1}\rtimes D_\chi\) by  $\widehat\theta(zd)=\theta(z)$ for any $z\in C_{q^2-1}$ and $d\in D_\chi$. On the other hand,
since $W$ is unipotent, Lemma \ref{lem:unipotent} shows that  $W$ has an $(\HH_p)_W$-invariant extension to ${}^2B_2(q^2)\rtimes D$. We write \(\widehat W\) for its restriction to \({}^2B_2(q^2)\rtimes D_\chi\). Define 
\[
\widehat\rho(zsd)=\theta(z)\widehat W(sd),
\]
where $z\in C_{q^2-1},\ s\in {}^2B_2(q^2)$ and $d\in D_\chi$.
Then  $\widehat\rho\in\Irr(M\rtimes D_\chi)$ and it is an extension of \(\theta\otimes W\). Inflate it via
$(J\rtimes D_\chi)/U\cong M\rtimes D_\chi$ to obtain \(\widehat\lambda\in\Irr(J\rtimes D_\chi)\).
Let $\widehat\chi=\Ind_{J\rtimes D_\chi}^{G\rtimes D_\chi}(\widehat\lambda)$. Then $\widehat\chi$ is irreducible.
Since \(G\rtimes D_\chi=G(J\rtimes D_\chi)\) and \(G\cap (J\rtimes D_\chi)=J\), the restriction yields
\[
\widehat\chi|_G=\Ind_J^G\Inf_M^J(\theta\otimes W)=\chi.
\]

Finally, we verify the stability of $\widehat \chi$. For any
$a=d\sigma\in(D\times\HH_p)_\chi$. Then by Equalities  \eqref{eq:hc} and \eqref{eq: act},  $\theta^{a}=\theta^{n}\in\{\theta,\theta^{-1}\}$
and $W^{a}=W^{n}=W$. Since \(W^d=W\) by Lemma \ref{defect zero characters of $K$}, we have \(W^\sigma=W\). Therefore, the chosen extension satisfies $\widehat W^{a}=\widehat W$. If $\theta^{a}=\theta$, then by $\widehat \rho$ immediately yields
$\widehat\lambda^{\,a}=\widehat\lambda$. The induction commutes with the action of $a$, so $\widehat\chi^{a}=\widehat\chi$. If
$\theta^{a}=\theta^n=\theta^{-1}$. Since \(n\) centralizes \({ }^2B_2(q^2)\) and is $D$-fixed, it also centralizes \({}^2B_2(q^2)\rtimes D_\chi\). Hence, on the group \(M\rtimes D_\chi\), we obtain
\[
\widehat\rho^{\,a}(zsd)=
\theta^{-1}(z)\widehat W(sd)=
\widehat\rho^{\,n}(zsd).
\]
Moreover, since $n^d=n$ and $J\cap J^n=M$, 
$(J\rtimes D_\chi)\cap(J\rtimes D_\chi)^{\,n}=
M\rtimes D_\chi$. Apply Mackey's formula and Frobenius reciprocity (see \cite[Section 5]{I1976}),
we obtain
\[
\begin{aligned}
\left\langle
\widehat\chi^{\,a},\widehat\chi
\right\rangle_{G\rtimes D_\chi} &=
\sum_{x\in J\rtimes D_\chi\backslash G\rtimes D_\chi/J\rtimes D_\chi}
\left\langle
\widehat\lambda^{\,a}|_{(J\rtimes D_\chi)\cap (J\rtimes D_\chi)^x},\widehat\lambda^{\,x}|_{(J\rtimes D_\chi)\cap (J\rtimes D_\chi)^x}
\right\rangle_{(J\rtimes D_\chi)\cap (J\rtimes D_\chi)^x},
\end{aligned}
\]
where the term corresponding to \(x=n\) equals
$\left\langle
\widehat\rho^{\,a},\widehat\rho^{\,n}
\right\rangle_{M\rtimes D_\chi}=
\left\langle
\widehat\rho^{\,n},\widehat\rho^{\,n}
\right\rangle_{M\rtimes D_\chi}=1$, and all other terms are non-negative integers. Thus 
$\left\langle
\widehat\chi^{\,a},\widehat\chi
\right\rangle_{G\rtimes D_\chi}\geq1$.
However, we know that $\widehat\chi^{\,a}$ and $\widehat\chi$
are irreducible, and this forces  $\widehat\chi^{\,a}=\widehat\chi$. We complete the proof.
\end{proof}

\begin{lemma}\label{lem: The non-split Suzuki-centralizer families}
If $\chi$ is any of characters in the families $\chi_{43},\chi_{44},\chi_{47},\chi_{48}$, then it has an $(\HH_p\times D)_\chi$-invariant extension  to $G\rtimes D_\chi$.
\end{lemma}
\begin{proof}
Take the torus $T_\pm\cong C_{q^2\pm\sqrt{2}q+1}\times C_{q^2\pm\sqrt{2}q+1}$. In the notation of \cite[Table~B.2]{HH}, these are $T_6$ and $T_7$, and the torus $C_{q^2+\sqrt{2}q+1}\times C_{q^2-\sqrt{2}q+1}$ is $T_5$ with order $q^4+1$. For a non-trivial $\theta\in\Irr(T_\pm)$, the formulas in \cite[Appendix~B, construction of the $g_8$ and $g_{10}$ series]{HH} give
\begin{align}
 R_{{\bf T}_-}^{\bf G}(\theta)&=-\chi_{42}+\chi_{43}+\chi_{44}+\chi_{45},\label{eq:dlminus}\\
 R_{{\bf T}_+}^{\bf G}(\theta)&=-\chi_{46}+\chi_{47}+\chi_{48}+\chi_{49}.\label{eq:dlplus}
\end{align}
For convenience, from left to right, we will denote these four constituents as $\chi_1,\chi_2,\chi_3,\chi_4$.

First, we claim that $T_\pm$ is stable under a representative element $F_0$ of $D$, and that all rational Weyl representative elements can be taken $D$-fixed.
Let $\gamma$ be the exceptional endomorphism with $\gamma^2=F_2$. In the twisted case, write $F'=c_n\gamma^{2m+1}$ and $v=c_n\bar\gamma$, where $\bar{\gamma}$ is a finite-order action on the root system after removing the Frobenius scaling and $c_n$ 
denotes the inner automorphism action of $n$. Here, $n$ is
a preimage of an element of the Weyl group and $v$ has order $8$. Let $k=2m+1$. Since \(k\) is odd, there exists an integer \(b\) such that $bk\equiv1\pmod8$.
Hence, $(v^b)^k=v^{bk}=v$. Also, since \(b\) is odd, we may write $v^b=n_0\bar\gamma$, and define 
$F_0=c_{n_0}\gamma$. Then in characteristic $2$, 
the Tits kernel $\langle h_\alpha(-1)\rangle$ is trivial, so
we deduce that $F_0^k=F'$ and $F_0$ fixes all chosen  rational Weyl representative elements. That is, $F_0$ fixes $W^{F'}=\CC_{W}(v)$ pointwise, where $W=\NN_G(T_{\pm})/T_\pm$ is the rational Weyl group of $G$, as claimed.

Next, for every $\chi$ lying in the families $\chi_i(i=1,2,3,4)$, we  will construct an extension $\widehat \chi$ of $\chi$ to $G\rtimes D_\chi$ through the following paragraphs.

Write $H=W_\theta$ for the stabilizer of $\theta$ in $W$. Note that $\theta$ is trivial on $Z_{p'}=(C_{ q^2\pm\sqrt{2}q+1})_{p'}$,  so by Lemma \ref{E}, we have $H\cong C_4$. Also, by Table \ref{tab:1}, $W\cong G_8$ which has order $96$. Then $\ZZ(W)\cong C_4$. By \cite[Theorem 5.5]{MT18}, $\NN_W(H)=H\times \ZZ(W)$, so $\NN_W(H)/H\cong C_4$. Let \(d\in D_\chi\). Then \(d\) fixes the corresponding $G$-class of a semisimple element, and so we may write $\theta^d=\theta^{w_d}$ for some $w_d\in W$.
Thus  $W_{\theta^d}=W_\theta=H$ and \(w_d\) normalizes \(H\). This yields a group homomorphism $D\rightarrow\NN_W(H)/H$, and the order of its image divides both \(|D_\chi|\) and \(4\). However, \(|D_\chi|\) is odd, so this image must be trivial. Consequently, for any $d\in D_\chi$, we have $\theta^d=\theta$. Define
$\widehat\theta(td)=\theta(t)$,
where $t\in T_{\pm}$ and $ d\in D_\chi$. Clearly, $\widehat \theta$ is a linear character of $T_{\pm}\rtimes D_\chi$ and  
an extension of $\theta$.

Since the action of $D$ fixes the chosen Weyl representative elements, every element of \(\NN_G(T_{\pm})\) can be written as $tn_w$, where $t\in T_{\pm}$ and $w\in W$ with each \(n_w\) commuting with \(D\). Thus
$T_{\pm}\rtimes D\unlhd \NN_G(T_\pm)\rtimes D$ 
and  $\NN_G(T_\pm)\rtimes D/T_{\pm}\rtimes D\cong W$.
By Frobenius reciprocity and Mackey's formula(see \cite[Section 5]{I1976})
\begin{equation}\label{2}
\begin{aligned}
\left\langle
\Ind_{T_\pm\rtimes D_\chi}^{\NN_G(T_\pm)\rtimes D_\chi}(\widehat \theta),
\Ind_{T_\pm\rtimes D_\chi}^{\NN_G(T_\pm)\rtimes D_\chi}(\widehat \theta)
\right\rangle_{\NN_G(T_\pm)\rtimes D_\chi}
&=\sum_{w\in W}
\left\langle
\widehat\theta,\widehat\theta^{\,n_w}
\right\rangle_{T_\pm\rtimes D_\chi}.
\end{aligned}
\end{equation}
Since \(\widehat\theta\) is linear, each summand is either \(0\) or \(1\). Furthermore, since \(n_w\) commutes with \(D\), and  both sides are trivial on \(D_\chi\), so $\widehat\theta^{\,n_w}=\widehat\theta$ if and only if  $\theta^w=\theta$. Thus Equality \eqref{2} is equal to $|H|=4$.

To construct the required extension, we need to realize $G\rtimes D_\chi$ as a disconnected fixed-point group, and form the disconnected Deligne-Lusztig virtual character $\widehat R$ from this extension $\widehat \theta$. Let $|D_\chi|=e$, and write $r=\frac{2m+1}{e}$ and $F_1=F_0^r$. Then $F_1^e=F'$. 
Write ${\bf X}={\bf G}^e$, and let $\tau$ cyclically permute these $e$ factors by $\tau(x_1,\dots,x_e)=(x_e,x_1,\dots,x_{e-1})$. Define 
\[
{\cal F}(x_1,\dots,x_e)=(F_1(x_e),F_1(x_1),\dots,F_1(x_{e-1}))
\]
Then $x\mapsto (x,F_1(x),\dots,F_1^{e-1}(x))$
gives  $G={\bf G}^{F'}\cong{\bf X}^{\cal F}$. Moreover, $\cal F$ is also a Steinberg map.   Under this correspondence, the action of $\tau$ corresponds to the action of $F_1^{-1}$. Therefore, write
${\bf Y}={\bf X}\rtimes\langle\tau\rangle$ and let ${\cal F}$ fix $\tau$. We then obtain the finite group structure
\begin{equation}\label{iso}
G\rtimes D_\chi\cong {\bf Y}^{\cal F}.    
\end{equation}
Take a suitable $F_1$-stable maximal torus ${\bf T}$ in ${\bf G}$ contained in the $F_1$-stable Borel subgroup ${\bf B}$. 
Then the disconnected torus ${\bf L}={\bf T}^e\rtimes \langle \tau\rangle$ is well-defined and $\cal F$-stable.
Thus $T_{\pm}\rtimes D_\chi\cong{\bf L}^{\cal F}$ by the isomorphism $T_{\pm}\cong({\bf T}^e)^{\cal F}$.
As in \cite[Definition 2.4]{DM}, we choose 
$\widehat t\in\Irr({\bf L}^{\cal F})$ corresponding to $\hat{\theta}\in\Irr(T_{\pm}\rtimes D_\chi)$, and then construct the disconnected Deligne-Lusztig virtual character
\[
\widehat R:=R_{{\bf L}}^{{\bf Y}}(\widehat t)\in{\mathbb Z}\Irr({\bf Y}^{\cal F}).
\]
By the isomorphism \eqref{iso}, 
note that $\widehat t$ is an extension of $t$ because 
$\widehat \theta$ is an extension of $\theta$.
Since the restriction commutes with the Deligne-Lusztig induction(see \cite[Proposition 3.2.9]{GM20}), it is not difficult to verify that
\begin{equation}\label{com}
{\rm Res}^{{\bf Y}^{\cal F}}_{{\bf X}^{\cal F}}(\widehat R)
=
R_{{\bf T}^e}^{\bf X}\big({\rm Res}_{({\bf T}^e)^{\cal F}}^{{\bf L}^{\cal F}}(\widehat t)\big)
=R_{{\bf T}^e}^{{{\bf X}}}(t)
=R_{{\bf T}_\pm}^{\bf G}(\theta)
=-\chi_1+\chi_2+\chi_3+\chi_4,
\end{equation}
where $t\in\Irr(({\bf T}^e)^{\cal F})$ corresponds to  $\theta\in\Irr(T_\pm)$.  Moreover, we obtain
\begin{equation}\label{4}
\begin{aligned}
\langle \widehat R, \widehat R \rangle_{{\bf Y}^{\cal F}}
=\left\langle R_{{\bf L}}^{{\bf Y}}(\widehat t),  R_{{\bf L}}^{{\bf Y}}(\widehat t)\right\rangle_{{\bf Y}^{\cal F}}
&=\left\langle \Ind_{{\bf L}^{\cal F}}^{\NN_{{\bf Y}^{\cal F}}({\bf T}^e)}(\widehat t), \Ind_{{\bf L}^{\cal F}}^{\NN_{{\bf Y}^{\cal F}}({\bf T}^e)}(\widehat t) \right\rangle_{\NN_{{\bf Y}^{\cal F}}({\bf T}^e)}\\
&=
\left\langle
\Ind_{T_\pm\rtimes D_\chi}^{\NN_G(T_\pm)\rtimes D_\chi}(\widehat \theta),
\Ind_{T_\pm\rtimes D_\chi}^{\NN_G(T_\pm)\rtimes D_\chi}(\widehat \theta)
\right\rangle_{\NN_G(T_\pm)\rtimes D_\chi}=4,
\end{aligned}
\end{equation}
where the second equality by \cite[Proposition 3.2]{DM}, the third equality induced by \eqref{iso} and the last equality by \eqref{2}.

For every $\chi_i$, which comes from Equality \eqref{eq:dlminus} or \eqref{eq:dlplus}. By the degrees in \cite[Table~A.2]{Himstedt},  we can partition them into three parts $\{\chi_1\}, \{\chi_2,\chi_3\}, \{\chi_4\}$, where the middle pair have the same degree, other parts have distinct degrees. So $\chi_i^d=\chi_i$ for any $d\in D$ and $i=1,4$. Also, since $D$ has odd order, it must act trivially on $\{\chi_2,\chi_3\}$. Thus every $\chi_i$ is $D$-fixed.
Since $D_\chi$ is cyclic, $\chi_i$ can be extended to $G\rtimes D_\chi$. Choose an extension $\chi_i^*$ and  for any  $\lambda\in\Irr(D_\chi)$, all $\chi_i^*\lambda$ are extensions of $\chi_i$. Therefore, together with Equality \eqref{com} we may write
\[
\widehat{R}=\sum_{i=1}^4\sum_{\lambda\in\Irr(D_\chi)}a_{i,\lambda}\chi_i^*\lambda+\Psi
\]
with $a_{i,\lambda}\in\mathbb{Z}$, and every irreducible constituent of $\Psi$ does not lie above any of the four $\chi_i$. Restricting to $G$ and comparing the coefficients of $\chi_i$, we obtain
$\sum_{\lambda}a_{i,\lambda}=\varepsilon_i=\pm 1$.
Note that for each $i$, $\sum_{\lambda}a_{i,\lambda}^2\ge 1$.
On the other hand, together with Equality \eqref{4}, we have
\[
4=\langle\widehat{R},\widehat{R}\rangle_{G\rtimes D_\chi}
=\sum_{i=1}^4\sum_{\lambda}a_{i,\lambda}^2+\langle \Psi,\Psi\rangle_{G\rtimes D_\chi}\ge 4.
\]
This forces $\Psi=0$ and $\sum_{\lambda}a_{i,\lambda}^2=1$ for every $i$. So $\widehat{R}=\sum_{i=1}^4\varepsilon_i\widehat{\chi}_i$,
$\widehat{\chi}_i\big|_G=\chi_i$ and every $\widehat{\chi}_i$ is the unique irreducible character in $\widehat{R}$ lying above $\chi_i$. 

Finally, we verify the stability. Let
$a=d\sigma\in (D\times \HH_p)_\chi$. 
Then $a$ fixes the corresponding Lusztig series.
So there exists some $w\in W$ such that $\theta^a=\theta^w$.
Since $w$ admits a $D$-fixed representative element $n_w\in\NN_G({\bf T})$, and $\widehat{\theta}$ is trivial on $D_\chi$, we have $\widehat{\theta}^a=\widehat{\theta}^{n_w}$.  Two compatibilities for Deligne-Lusztig virtual character ${\widehat R}$ on the above disconnected group $({\bf Y}, \cal F)$ are required
\begin{equation}\label{7}
R_{{\bf L}}^{{\bf Y}}(\widehat t)^a=R_{{\bf L}}^{{\bf Y}}(\widehat{t}^a)
\quad \text{and}\quad R_{{\bf L}}^{{\bf Y}}(\widehat{t}^{n})=R_{{\bf L}}^{{\bf Y}}(\widehat{t}),
\end{equation}
where $n=(n_w,\ldots,n_w)\in{\bf X}^{\cal F}$ 
corresponds to $n_w$. By definition of ${R}_{\bf L}^{\bf Y}$, and together with \cite[Proposition 3.3]{DL}, the first equality holds. Since $n_w\in\NN_G({\bf T})$ and $F_0(n_w)=n_w$,
we have  $n{\bf T}^en^{-1}={\bf T}^e$ and $n\tau=\tau n$, so that $n\mathbf L n^{-1}=\mathbf L$. 
Consider the inner automorphism $c_n:\mathbf Y\rightarrow\mathbf Y, x\mapsto nxn^{-1}$.
Note that \(n\) commutes with \(\cal F\) and induces an isomorphism of the corresponding algebraic varieties.
Using the independence of the quasi-parabolic choice(see \cite[Theorem 4.8.13]{GM20}), and since $R_{{\bf L}}^{{\bf Y}}(\widehat t)$ is a class function, 
$R_{{\bf L}}^{{\bf Y}}(\widehat{t}^{n})=R_{{\bf L}}^{{\bf Y}}(\widehat{t})^n=R_{{\bf L}}^{{\bf Y}}(\widehat{t})$.
Return to the group $G\rtimes D_\chi$, by Equality \eqref{7} and
$\widehat{\theta}^a=\widehat{\theta}^{n_w}$, we obtain
\begin{equation}\label{8}
R_{T_\pm\rtimes D_\chi}^{G\rtimes D_\chi}(\widehat \theta)^a=R_{T_\pm\rtimes D_\chi}^{G\rtimes D_\chi}(\widehat \theta^a)=R_{T_\pm\rtimes D_\chi}^{G\rtimes D_\chi}(\widehat \theta^{n_w})=R_{T_\pm\rtimes D_\chi}^{G\rtimes D_\chi}(\widehat \theta).    
\end{equation}
Also, since \(\widehat\chi\) is the unique irreducible character in \(\widehat R\) lying over \(\chi\) and
$(\widehat\chi^a)|_G=(\widehat\chi|_G)^a=\chi^a=\chi$.
Meanwhile, Equality \eqref{8} guarantees that \(\widehat\chi^a\) still appears in \(\widehat R\) with the same coefficient. Therefore, the uniqueness forces $\widehat\chi^a=\widehat\chi$. We complete the proof.
\end{proof}

All extensions of the weight characters are handled by the following proposition. 
\begin{proposition}\label{equ:local extention p>=5}
Let $\chi\in\IBr_p(G)$ and  $\Omega(\chi)=\overline{(R,\varphi)}\in\WW_p(G)$.
There exists $D_0\leq D$ such that $\varphi^\circ\in\IBr(\NN_G(R))$  has an $(\HH_p\times D_0 )_\varphi$-invariant extension to $\NN_{G}(R)\rtimes (D_0)_{\varphi}$. 
\end{proposition}
\begin{proof}
Suppose that $R$ is non-cyclic abelian. Then by Table \ref{tab:1}, $\NN_G(R)/R=T_{p'}\rtimes W$, where $W$ is isomorphic to  $D_{16}$, $\GL_2(3)$, or $G_8$. Note that $p\geq 5$  and $p\nmid |\NN_G(R)/R|$, so $p\nmid\varphi(1)$. Thus we use \cite[Proposition~4.5]{Johansson},
$\varphi$ has an $(\HH_p\times \NN_{G\rtimes D}(R))_\varphi$-invariant extension $\widehat\varphi\in\Irr(\NN_{G\rtimes D}(R)_\varphi)$. The restriction $\widehat \varphi^\circ$ is
an $(\HH_p\times \NN_{G\rtimes D}(R))_\varphi$-invariant extension of $\varphi^\circ$.

Suppose that $R$ is cyclic. For convenience, we write $N=(Z_{p'}\rtimes E)\times K$ and $I=(Z_{p'}\rtimes E_\theta)\times K$. Then Equality \eqref{eq:rankone} gives 
\[
\varphi:=\psi_{\theta,\kappa,\epsilon}=
 \Ind_{I}^{N}
((\widehat\theta\epsilon)\otimes\kappa)\in\dz_p(\NN_G(R)/R),
\]
where $\widehat \theta\in\Irr(Z_{p'}\rtimes E_\theta)$ is
the trivial extension of $\theta\in\Irr(Z_{p'})$, $\epsilon\in\Irr(E_\theta)$ and $\kappa\in\dz_p(K)$.
Moreover, in the case, $\NN_G(R)=(Z\rtimes E)\times K$ and
$[E,K]=1$, where $Z$ and $K$ are as in Table \ref{tab:2}, and $E=C_2$ or $C_4$.  By Clifford Theorem, we obtain $\operatorname{Res}_{Z_{p'}\times K}^{N}(\varphi)=\sum_{e\in E/E_\theta}\theta^e\otimes\kappa$.

Let $D_0\leq D$ be as in Lemma \ref{lem:3.5}. That is, $D_0$ has odd order, $D_0$ normalizes $Z_{p'}$ and $K$, and centralizes $E$. For any $d\in (D_0)_\varphi$, $\kappa^d=\kappa$ and $\theta^d=\theta^e$ for some $e\in E$. 
Since $[d,E]=1$, $d$ commutes with the $E$-action on $\Irr(Z_{p'})$. Consequently, the map
$d\mapsto eE_\theta$ defines a group homomorphism $(D_0)_\varphi\rightarrow E/E_\theta$. However, $(D_0)_\varphi$
has odd order and $E/E_\theta$ is $2$-group. So this homomorphism must be trivial. Thus we also have $\theta^d=\theta$.

First, we construct an $(\HH_p\times D_0)_\kappa$-invariant extension $\widehat\kappa$ of $\kappa$ to  $K\rtimes (D_0)_\varphi$ as follows.  By Table \ref{tab:2}, $K=\SL_2(q^2)$ or ${}^2B_2(q^2)$. For the former, all $\kappa$ are real-valued because $q^2$ is even.  For the latter, by \cite[p.71 and Table 4.2.1]{S16}, all $\kappa$ are real-valued expect for two unipotent characters $W^{+}$ and $W^{-}$ when $p\mid q^2-1$. Suppose that $\kappa$ is real-valued. Note that $(K\rtimes (D_0)_\varphi)/K\simeq (D_0)_\varphi$ is cyclic with odd order. Then by Lemma \ref{lem:real}, $\kappa$ has the unique real-valued extension $\widehat{\kappa}\in\Irr(K\rtimes (D_0)_\varphi)$. For any $a=\sigma d\in(\HH_p\times D_0)_\kappa$, \(\widehat\kappa^{a}\) is still a real-valued extension of \(\kappa\). Thus the uniqueness yields
$\widehat\kappa^{a}=\widehat\kappa$. Suppose that $\kappa=W^{\pm}$. In this case, by Lemma \ref{lem:unipotent}, we can choose $\widetilde\kappa\in\Irr(K\rtimes D_0)$ such that it is an \((\HH_p)_\kappa\)-invariant extension of \(\kappa\). Let $\widehat\kappa=\operatorname{Res}^{K\rtimes D_0}_{K\rtimes  (D_0)_\varphi}(\widetilde\kappa)$. Then it is also an extension of \(\kappa\). Since $\kappa^d=\kappa$,
$\kappa^{a}=\kappa$ implies \(\kappa^\sigma=\kappa\).
On the other hand, the $d$-conjugation in \(K\rtimes D_0\) is an inner automorphism, and thus $d$ fixes \(\widetilde\kappa\). Consequently, we also have $\widehat\kappa^{a}=\widehat\kappa$. 

Next, we construct an extension $\widehat \varphi$ of $\varphi$ to $N\rtimes (D_0)_\varphi$. Write $\widehat I=I\rtimes (D_0)_\varphi$ and $\widehat N=N\rtimes (D_0)_\varphi$. Since $\theta$ is $(D_0)_\varphi$-fixed and $(D_0)_\varphi$ centralizes $E$, we may write $\widehat I=(Z_{p'}\times K)\rtimes (E_\theta \times (D_0)_\varphi)$. Define  $\widetilde{\theta}(ted)=\theta(t)$ and 
$\widehat{\epsilon}(ed)=\epsilon(e)$, where 
$t\in Z_{p'}$, $e\in E_\theta$ and $d\in (D_0)_\varphi$.
Then $\widetilde \theta$ and $\widehat \epsilon$ are linear characters of $Z_{p'}\rtimes(E_\theta\times (D_0)_\varphi)$ and  $E_\theta\times (D_0)_\varphi$, respectively.
Define on $\widehat I$ 
\[
\widehat\beta(tked)=\theta(t)\epsilon(e)\widehat\kappa(kd),
\]
where $t\in Z_{p'}$, $k\in K$, $e\in E_\theta$, $d\in (D_0)_\varphi$.
Then  $\widehat\beta$ is a character of $\widehat I$.
Set $\beta=\widehat \theta\epsilon\otimes \kappa$,
and then $\operatorname{Res}_{I}^
{\widehat I}(\widehat\beta)=\beta$. Since $\beta$ is irreducible, $\widehat\beta$ is also irreducible. Now write
$\widehat\varphi=\Ind_{\widehat I}^{\widehat N}(\widehat\beta)$.
Then $\widehat\varphi$ is a character of $\widehat N$.
Since $\widehat N=\widehat IN$ and $N\cap \widehat I=I$, 
the Mackey's formula(see \cite[Problem (5.2)]{I1976}) gives
\[
\operatorname{Res}_{N}^{\widehat N}(\widehat\varphi)
=\operatorname{Res}_{N}^{\widehat N}(\Ind_{\widehat I}^{\widehat N}(\widehat\beta)) 
=\Ind_{I}^{N}(\operatorname{Res}_{I}^
{\widehat I}(\widehat\beta))
=\Ind_{I}^{N}(\beta)=\varphi.
\]
Since $\varphi$ is irreducible, so $\widehat \varphi$  is indeed an extension of $\varphi$.

Finally, we verify the stability of $\widehat \varphi$.
Let $a=d\sigma\in(D_0\times\HH_p)_\varphi$. The direct product decomposition and the uniqueness of Clifford correspondence, there exists some $e_a\in E$ such that $\theta^a=\theta^{e_a}$, $\kappa^a=\kappa$ and $\epsilon^\sigma=\epsilon$.
Here, the last equality can also be obtained by Lemma \ref{lem:3.5}. From the earlier choice of \(\widehat\kappa\), we have $\widehat\kappa^a=\widehat\kappa$.
Since $[E,(D_0)_\varphi]=1$ and $ [E,K]=1$, \(e_a\) normalizes \(I\) and acts trivially on the factor \(K\rtimes (D_0)_\varphi\). Therefore, by the definition of $\widehat \beta$, we obtain $\widehat\beta^a=\widehat\beta^{\,e_a}$.
Compatibility of induction with the action of $a$, this gives 
\[
\widehat\varphi^a
=\Ind_{\widehat I}^{\widehat N}(\widehat\beta^a)
=\Ind_{\widehat I}^{\widehat N}(\widehat\beta^{e_a})
=(\Ind_{\widehat I}^{\widehat N}(\widehat\beta))^{e_a}
=\widehat\varphi^{e_a}
=\widehat \varphi,
\]
where the final equality holds because $e_a\in N \leq {\widehat N}$ and $\widehat \varphi$ is a class function on $\widehat N$. The restriction $\widehat \varphi^\circ$ gives the desired extension of $\varphi^\circ$.  We complete the proof.
\end{proof}

\subsection{Equivalent bijection and compatible extensions for  \texorpdfstring{$p=3$}{p=3}}\label{p=3}
In this subsection, keep the groups $G$ and $D$ as above and we fix $p=3$. Note that the prime $3$ is a bad prime for ${\bf G}$, and it only divides $q^2+1$ and $q^4-q^2+1$. Moreover, the Sylow $3$-subgroups of $G$ are non-abelian.

Let $q^2+1=3^a b$, where $(3,b)=1$ and $a\geq 1$. That is, $a=v_3(q^2+1)$ and $b=(q^2+1)_{3'}$. We will use the block and radical-normalizer classification in \cite[Tables~3--4]{Malle}, and together with the decomposition matrix in \cite[Section 4]{Himstedt}.
\subsubsection{The principal \texorpdfstring{$3$}{3}-block}
Let $B_0$ be the principal block of $G$. Then it follows from  \cite[first row of Table~3]{Malle} that $B_0$ has nine irreducible Brauer characters. The following notation is the CHEVIE notation. By \cite[Section 4 and Table~C.5]{Himstedt}, we write 
\begin{equation}
\IBr_3(B_0)=\{\phi_1,\phi_4,\phi_{5,1},\phi_7,\phi_8,\phi_{10},\phi_{15},\phi_{18},\phi_{21}\},
\end{equation}
which corresponds to the ordinary character set, forms a unitriangular set of $\Irr(B_0)$, as follows.
\begin{equation}\label{principal 3-block}
\mathcal{B}_0=\{\chi_1,\chi_4,\chi_{5,1},\chi_7,\chi_8,\chi_{10},\chi_{15},\chi_{18},\chi_{21}\}.
\end{equation}

\begin{lemma}\label{eq:nine}
All Brauer characters in block $B_0$ are rational-valued and $D$-invariant.
\end{lemma}
\begin{proof}
For $i\in\{1,4,7,8,10,18,21\}$ in \eqref{principal 3-block}, every \(\chi_i\) is unipotent by \cite[Section 4]{Himstedt}. Moreover, it follows from \cite[Table A.1]{Himstedt} that the degree of $\chi_i$ uniquely determines it. On the other hand, Galois automorphisms preserve degrees of characters. These 
characters \(\chi_i \) are rational-valued. 
For the non-unipotent character $\chi_{5,1}$(see \cite[Table A.2]{Himstedt}), by \cite[Table~C.5]{Himstedt}, it's easy to 
check that 
\[
\phi_{5,1}
=\chi_9^\circ-\phi_1+\phi_4+\phi_8-\phi_{10}.
\]
Note that by \cite[Table A.1]{Himstedt}, all terms on the right-hand side are rational-valued, so \(\phi_{5,1}\) is also rational-valued. For $\chi_{15}$, by \cite[Table~A.1]{Himstedt}, we have
$\QQ (\chi_{15})=\QQ(\sqrt{3})$. So $\chi_{15}$ is not rational-valued. However, $\chi_{15}^{\circ}=\chi_{16}^\circ$ is rational-valued. 
To see this, we take any $3'$-element \(g\in G\). Then $\chi_{15}(g)\in\QQ(\zeta_{|g|})$, and thus $\chi_{15}(g)\in
\mathbb Q(\sqrt{3})\cap\mathbb Q(\zeta_{|g|})=\mathbb Q$. 

Therefore, the unitriangularity implies that the nine Brauer characters are rational-valued. It follows from that \cite[Proof of Theorem 5.1]{Malle} that they are also fixed by every field automorphism.
\end{proof}

It follows from  \cite[Theorem 5.1]{Malle} that there are nine $3$-weights corresponding to nine $3$-Brauers in the block $B_0$ for $a=1$ or $a\geq 2$. As in \cite[Table 4]{Malle} and \cite[(2B)]{An}, non-trivial $3$-radical subgroups, corresponding normalizers and  $3$-defect  zero characters in Table \ref{tab:3}. Here,  $P_3$ denotes the Sylow $3$-subgroup of $G$, $Q_8$ is a quaternion group of order $8$, and ${-}.{-}$ denotes the extension of two groups. The last row because $b = (q^2 + 1)_{3'}$.
\begin{table}[ht]
\centering
\caption{Non-trivial \(3\)-radical subgroups and  $3$-defect  zero characters}\label{tab:3}
\begin{tabular}{@{}llllrr@{}}
\toprule
$R$ & $\NN_G(R)/R$ & The degrees of $\dz_3(\NN_G(R)/R)$& condition\\
\midrule
$P_3$ & $C_2\times C_2$ & $1,1,1,1$& $a\geq 2$\\
$3^{1+2}_+$ & $\SL_2(3)$ & $3$& $a\geq 2$\\
$3^{1+2}_+$ & $\SL_2(3)\rtimes C_2$ & $3,3$& $a\geq 2$\\
$3^{1+2}_+$ & $Q_8\rtimes C_2$ & $ 1,1,1,1,2,2,2$& $a=1$\\
$C_{3^a}\times C_{3^a}$ &$(C_b\times C_b)\rtimes \GL_2(3)$& $3,3$& $a\geq 1$ \\
\bottomrule
\end{tabular}
\end{table}
\begin{lemma}\label{lem:3.15}
Let $R$ be a radical subgroup in  Table \ref{tab:3}.
Then all defect zero characters of $\NN_G(R)/R$ are
$(\HH_3\times D)$-invariant.
\end{lemma}
\begin{proof}
Suppose that  $a\geq2$. For $C_2\times C_2$, all values of its four linear characters are integer. For $\SL_2(3)$, its the unique degree-$3$ character is inflated from $A_4\cong \PSL_2(3)$. So values of the degree-$3$ character are integer. For the last row when $a>1$, the two degree-$3$ characters are trivial on the normal subgroup \(C_b\times C_b\), and then they are inflated from the quotient group \(\mathrm{GL}_2(3)\).
Also, two degree-$3$ characters of $ \GL_2(3)$ are also  inflated from $S_4$. For \(H=\SL_2(3)\rtimes C_2\), write $K=\mathrm{SL}_2(3)$. Then $H/\ZZ(K)\cong S_4$. Thus two degree-$3$ characters of $H$ are inflated from $S_4$. Note that all values of characters of $S_4$ are integer. In short,
these characters are rational-valued.

Suppose that $a=1$. For the last row, similar to the above paragraph. These characters are rational-valued. For $Q_8\rtimes C_2$, seven characters come from 
\[
SD_{16}=\langle x,y\mid x^8=y^2=1,\ yxy^{-1}=x^3\rangle,
\]
and its every element is uniquely written as \(x^i\) or \(x^iy\) for \(i\in{\mathbb Z}/8{\mathbb Z}\).
For any $\sigma=\sigma_3^r\in\HH_3$ and $\theta\in \Irr(SD_{16})$, we have $\theta(g)^\sigma=\theta(g^{3^r})$ for $g\in\{x^i,x^iy\}$. However, note that $3^r\equiv1\ \text{or}\ 3\pmod8$, and one verifies directly that that every element of $SD_{16}$ is conjugate to its cube. So $\theta(g)^\sigma=\theta(g)$ for all $g$.
Therefore,  every character in $\Irr(SD_{16})$
is $\HH_3$-invariant.

Finally, it follows from that \cite[Proof of Theorem 5.1]{Malle} that $R$ is normalized by field automorphisms, and these field automorphisms stabilize the irreducible characters in Table \ref{tab:3}. Thus these characters are also $D$-invariant. We complete the proof.
\end{proof}

\subsubsection{The non-principal positive-defect \texorpdfstring{$3$}{3}-blocks}
Let \(s\in G\) be a non-trivial semisimple \(3'\)-element
and $B$ a non-principal positive-defect $3$-block of $G$.
Every character in \(\Irr_3(B)\) admits a Jordan decomposition parametrized by $s$ and a unipotent character $\nu$ of \(\CC_{G}(s)\). Write \(\chi_{s,\nu}\in\Irr(B)\) for the character with Jordan decomposition \((s,\nu)\), and let \(\eta_{s,\nu}\in\IBr_3(B)\) denote the $3$-Brauer character corresponding to \(\chi_{s,\nu}\) via  the chosen triangular basic set in \cite[Section 4]{Himstedt}.

It follows from \cite[Table 3]{Malle} that  $\CC_G(s)$ is a torus, or a direct product of $\SL_2(q^2)$ and a torus. Specifically, from the non-principal positive-defect $3$-blocks $2$ to $6$, they are respectively 
\begin{equation}\label{Cs}
C_{q^2+1}\times\SL_2(q^2),\
C_{q^2+1}\times C_{q^2+1},\
C_{q^2-1}\times\SL_2(q^2),\
C_{q^2-1}\times C_{q^2+1},\
C_{q^4-q^2+1}    
\end{equation}

\begin{lemma}\label{lem:3.19}
Let $B$ be a non-principal positive-defect $3$-block of $G$. Then every  character $\eta_{s,\nu}\in\IBr_3(B)$ is real-valued. Moreover, for any $\sigma\alpha\in\HH_3\times D$ with $\sigma$ mapping $3'$-th root of unity $\zeta$ to $\zeta^{3^r}$,
we have $\eta_{s,\nu}^{\sigma\alpha}=\eta_{\alpha(s^{3^r}),\nu}$.
\end{lemma}
\begin{proof}
If the centralizer \(\CC_{G}(s)\) is a tours as above \eqref{Cs}, then \(\IBr_3(B)\) has exactly one element, denoted by \(\eta_{s,1}\).
If \(\CC_{G}(s)\) has the factor \(\SL_2(q^2)\) as above, then the  corresponding two Brauer characters are $\eta_{s,1}$ and $\eta_{s,\St}$.
Note that $\chi_{s,v}$ corresponds to $\eta_{s,v}$, where 
$\nu\in\{1,\St\}$. For the generators $\sigma_3\in \HH_3$ and $\alpha\in D$, since $v$ is $\sigma_3\alpha$-fixed,  by \cite[Theorem  B]{Johansson} and \cite[Theorem 4.1]{Malle}, we obtain $\eta_{s,\nu}^{\sigma_3\alpha}=\eta_{\alpha^{-1}(s^{3}),\nu}$. 

Take an $F$-stable maximal torus $\bf T$ of $\bf G$ containing $s$. In the Weyl group of ${}^2F_4$, the longest element satisfies $w_0=-1$, so it acts on the torus by inversion
$t\mapsto t^{-1}$ for any $t\in {\bf T}$. Since $-1$ is a central element and fixed by the action defining the Ree group, it also lies in the rational Weyl group $W=\NN_G(T)/T$ of this torus $T={\bf T}^F$. Using Lang's theorem, there exists  $n\in \NN_{G}(\mathbf T)$ such that $nsn^{-1}=s^{-1}$. So
$(s^{-1},\overline{\nu})$ is $n$-conjugate to  $(s,\nu)$.
Thus $\overline{\chi_{s,\nu}}=\chi_{s^{-1},\bar \nu}=\chi_{s,\nu}$.
\end{proof}

By \cite[Tables 3-4]{Malle}, non-trivial $3$-radical subgroups and associated normalizers for all non-principal positive-defect $3$-blocks $B$ in Table \ref{tab:4} below.  The corresponding semisimple centralizers $\CC_G(s)$ are as in \eqref{Cs}. Here, $q^2+1=3^ab$ for $a\geq1$ and $(3,b)=1$. We will separately consider these cases.
\begin{table}[ht]
\centering
\caption{Non-trivial \(3\)-radical subgroups and normalizers}\label{tab:4}
\begin{tabular}{@{}lllrr@{}}
\toprule
$B$ &$R$ & $\NN_G(R)/R$ &condition \\
\midrule
2,3&$C_{3^a}\times C_{3^a}$ &$(C_b\times C_b)\rtimes \GL_2(3)$ & $a\geq 1$\\
4,5&$C_{3^a}$ & $\SL_2(q^2)\times (C_b\rtimes C_2)$ & $a\geq 2$\\
4,5&$C_3$ & $\PSU_3(q^2)\rtimes C_2$& $a=1$\\
6&$C_3$ & $\PSU_3(q^2)\rtimes C_2$& $a\geq1$\\
\bottomrule
\end{tabular}
\end{table}

We next consider the action of $\HH_3\rtimes D$ on the set $\dz_3(\NN_G(R)/R)$. Then we obtain
\begin{lemma}\label{lem:3.20}
Let $R$ be a radical subgroup of $G$ in  Table \ref{tab:4}.
Then all characters of $\dz_3(\NN_G(R)/R)$ are real-valued.
Moreover, the action of $\HH_3\rtimes D$ is compatible with the parameterization of these characters.
\end{lemma}
\begin{proof}
First,  for $R\cong C_{3^a}\times C_{3^a}$, let $T=C_{3^ab}\times C_{3^ab}$. Then it is a maximal torus of $G$ containing $R$. Moreover, $\NN_G(T)/T\cong \GL_2(3)$, $T_{3}=R$ and $T_{3'}=C_{b}\times C_{b}$.  Write $E=\GL_2(3)$ and fix the standard isomorphism $\Irr(T)\cong T$. Let $\theta_s\in\Irr(T_{3'})$ correspond to the non-trivial semisimple $3'$-element $s\in T$. Then $\theta_s\neq 1$.
Since $\CC_G(s)$ is $C_{3^ab}\times\SL_2(q^2)$ or $C_{3^ab}\times C_{3^ab}$ by Table \ref{tab:4} and \eqref{Cs},
it is easily seen that $E_{\theta_s}\cong \NN_{\CC_G(s)}(T)/T\cong C_2$ or $1$.

Define the extension $\widehat{\theta_s}\in\Irr(T_{3'}\rtimes E_{\theta_s})$ of $\theta_s$ by $\widehat{\theta_s}(te)=\theta_s(t)$ for any $t\in T_{3'}$ and $e\in E_{\theta_s}$. Let \(\varepsilon_s\) be the non-trivial linear character of \(E_{\theta_s}\cong C_2\). Define 
\begin{equation}\label{first}
\psi_{s}^+=\Ind_{T_{3'}\rtimes E_{\theta_s}}^{T_{3'}\rtimes E}(\widehat{\theta_s}),\quad \psi_{s}^-=\Ind_{T_{3'}\rtimes E_{\theta_s}}^{T_{3'}\rtimes E}(\widehat{\theta_s}\varepsilon_s),\quad \psi_{s}=\Ind_{T_{3'}}^{T_{3'}\rtimes E}(\theta_s),
\end{equation}
where $\psi_s$ comes from the case of \(E_{\theta_s}=1\). Clearly, all $3$-defect zero characters of $\NN_G(R)/R$ are given in these ways by \cite[Theorem 5.1]{Malle}. Note that the unique non-trivial element \(z=-I_2\in\ZZ(E)\) acts on  \(T_{3'}\) by inversion
\(s^z =s^{-1}\), so we obtain $\theta_{s^z}=\theta_{s^{-1}}=\theta_s^{-1}$.
Then $\widehat{\theta_{s^z}}=\widehat{\theta_{s^{-1}}}=\widehat {\theta_s}^{-1}$, and thus 
$\psi_{s}^\pm=\psi_{s^{-1}}^\pm=\overline{\psi_{s}^\pm}$ 
and $\psi_{s}=\psi_{s^{-1}}=\overline{\psi_{s}}$.
For $\sigma_3\in\HH_3$, note that $\widehat\theta_s^{\sigma_3}=\widehat{\theta_s^{\sigma_3}}=\widehat{\theta_{s^3}}$ and $\theta_s^{\sigma_3}=\theta_{s^3}$.  So $(\psi_{s}^\pm)^{\sigma_3}=\psi_{s^{3}}^\pm$ 
and $\psi_{s}^{\sigma_3}=\psi_{s^{3}}$. 

Second, for $R\cong C_{3^a}$ and $a\geq 2$, by Table \ref{tab:4}, $\NN_G(R)/R\cong \SL_2(q^2)\times (C_{b}\rtimes C_2)$. Since $3\nmid 2b$, every $3$-defect zero character
of $\NN_G(R)/R$ is denoted as
\begin{equation}\label{3.19}
\psi_{\kappa,\beta}=\kappa\otimes\beta,
\end{equation}
where $\kappa\in\dz_3(\SL_2(q^2))$ and $\beta\in\Irr(C_{b}\rtimes C_2)$.
On the hand, note that \(|\SL_2(q^2)|=q^2(q^2+1)(q^2-1)\),
where $q^2$ is even and $v_3(q^2+1)=a$. 
So $v_3(|\SL_2(q^2)|)=v_3(q^4-1)=a$, and then the degree-\((q^2+1)\) characters are $3$-defect zero. Thus
$\kappa$ can be induced from  a non-trivial character $\lambda$ of the non-split torus $C_{q^2+1}$. That is, $\kappa=R_{C_{q^2+1}}^{\SL_2(q^2)}(\lambda)$, where $R$ denotes the Deligne-Lusztig induction. We write $\kappa=\kappa_{\lambda}$. Clearly, by the involution $w$ of the Weyl group, we have $\lambda^{w}=\lambda^{-1}$. Then $\kappa_{\lambda}=\kappa_{\lambda^{-1}}=\overline{\kappa_\lambda}$. Also, by \cite[Lemma 3.1]{S20}, we have $\kappa^{\sigma_3}=\kappa_{\lambda}^{\sigma_3}=\kappa_{\lambda^{\sigma_3}}$ for $\sigma_3\in\HH_3$.

On the other hand, since \(b\) is odd, $\Irr(C_{b}\rtimes C_2)$ consists of two linear characters $1$, $\varepsilon$, lifted from $C_2$, and $(b-1)/2$ characters $\rho_j$ of degree \(2\), indexed by the orbits \(\{j,-j\}\) of non-zero elements \(j\in\mathbb Z/b\mathbb Z\). We write $j\sim-j$. When $\beta\in\{1,\varepsilon\}$, $\bar\beta=\beta$ and $\beta^{\sigma_3}=\beta$ for $\sigma_3\in\HH_3$. Otherwise, $\beta$ can be induced from a character $\rho_j$ of $C_{b}$, where $\rho_j(x)=\zeta_b^j$. However, since $j\sim{-j}$, we have $\rho_j=\rho_{-j}=\overline{\rho_j}$.  Clearly, $\rho_j^{\sigma_3}=\rho_{3j}$.
Therefore, 
$\overline{\psi_{\kappa,\beta}}=\psi_{\bar{\kappa},\bar{\beta}}=\psi_{\kappa,\beta}$ and $\psi_{\kappa,\beta}^{\sigma_3}=\psi_{\kappa^{\sigma_3},\beta^{\sigma_3}}$.

Next, for $R\cong C_3$, let $S=\PSU_3(q^2)$ and $\tau$ is a graph automorphism of order $2$. Note that the order $|S|=q^6(q^6+1)(q^4-1)/3$. Moreover,  $v_3(q^2-1)=0$ and $v_3(q^4-q^2+1)=1$. Also, $q^4-1=(q^2+1)(q^2-1)$ and $q^6+1=(q^2+1)(q^4-q^2+1)$, so $v_3(q^4-1)=a$ and $v_3(q^6+1)=a+1$. Thus $v_3(|U|)=-1+(a+1)+a=2a$. This implies that $\theta\in\dz_3(S)$ if and only if $v_3(\theta(1))=2a$.

As in \cite[Sections~2.1 and~3.1]{Orevkov}, $8$ families of ordinary character degrees of ${\rm GU}_3(q^2)$ are listed,
using their notations, as follows.
\[
\{1,q^2(q^2-1),q^6,q^4-q^2+1,\ q^2(q^4-q^2+1),\ (q^2-1)(q^4-q^2+1),\ d_7=q^6+1,\ d_8=(q^2-1)(q^2+1)^2\},
\]
and after restricting to \({\rm SU}_3(q^2)\), other two degrees generated $(q^2-1)(q^4-q^2+1)/3$ and $d_8/3$.
It is easily checked that all members of $\dz_3(S)$ come from the two families of degrees $d_7$ and $d_8$. For $d_7$-family,  $v_3(d_7)=a+1$. For $d_8$-family, $v_3(d_8)=2a$ and $v_3(d_8/3)=2a-1$. The other families, including the restriction generated degrees(see \cite[Section 3.1]{Orevkov}), none of them has a \(3\)-defect zero character. 
Therefore, when \(a=1\), all $3$-defect zero characters of $S$ come exactly from the two families \(d_7,d_8\). When \(a\ge 2\),
all $3$-defect zero characters of $S$ come only from the \(d_8\) family. 

Suppose that $a=1$. We will give the parametrization  of $\dz_3(S)$ for the $d_7$-family as in \eqref{3.20}.
To avoid confusion, we write \(d_7\)-family in \cite{Orevkov} as $X_{t,s}$, where $t\in{\mathbb Z}/(q^2+1)$, $s\in{\mathbb Z}/(q^4-1)$, $s\not\equiv0\pmod{q^2-1}$ and 
$X_{t,s}=X_{t,-q^2s}$. Multiplying by the linear character $\lambda$ obtained from the determinant changes the parametrization via
$(t,s)\longmapsto(t+1,s+1-q^2)$. After restriction to \({\rm SU}_3(q^2)\), $\lambda=1$. Thus we may normalize $s$ to \(0\) and retain  $k=s+(q^2-1)t\pmod{q^4-1}$. The equivalence relation after restriction becomes $k\sim -q^2k$, and thus
$X_{0,k}$ denotes characters of the \(d_7\)-family.

We now check whether we can descend to \(S\). 
Let $eI_3\in \ZZ({\rm SU}_3(q^2))$ be a generator. 
Then it acts  on the corresponding representation by \(e^k\). 
So $\ZZ({\rm SU}_3(q^2))\leq\ker(X_{0,k}|_{{\rm SU}_3(q^2)})$
if and only if $3\mid k$. We write $k=3j$, where 
$j\in{\mathbb Z}/\frac{q^4-1}{3}=\mathbb Z/((q^2-1)b)$.
Since $(q^2-1,b)=1$, the Chinese Remainder Theorem gives a bijection
\[ 
j\longleftrightarrow (u,z)\in {\mathbb Z}/(q^2-1)\times{\mathbb Z}/b.
\]
The condition \(k\not\equiv0\pmod{q^2-1}\) is equivalent to \(u\neq0\). Since $k\sim -q^2k$, $-q^2\equiv-1\pmod{q^2-1}$ and $-q^2\equiv1\pmod b$, we obtain $(u,z)\sim(-u,z)$. Therefore,
we define 
\begin{equation}\label{3.20}
\{\rho_{u,z}|\ 0\ne u\in\mathbb Z/(q^2-1),\ z\in\mathbb Z/b,\ (u,z)\sim(-u,z)\}.    
\end{equation}
These characters are $3$-defect zero characters of $S$ for the $d_7$ family.

For $\tau$, which can be realized via inverse transposition. 
So $\rho_k^\tau=\rho_{-k}=\overline{\rho_k}$. Note that
$-k\longleftrightarrow(-u,-z)\sim(u,-z)$, so
$\rho_{u,z}^{\tau}=\rho_{u,-z}=\overline{\rho_{u,z}}$.
Moreover,  $\rho_{u,z}^{\sigma_3}=\rho_{3u,3z}$.
If $z=0$, then $\rho_{u,0}^{\tau}=\rho_{u,0}=\overline{\rho_{u,0}}$. 
So two defect zero characters of $N$, come from two extensions of $\rho_{u,0}$ to $S\rtimes \langle \tau \rangle$, are also real-valued.  Also, note that the determinants of these two extensions are $1$ and $\varepsilon$, the latter is the non-trivial linear character of $U\rtimes\langle\tau\rangle/U\cong C_2$. We write  $\psi_{u,0}^+, \psi_{u,0}^-$ for  two extensions of $\rho_{u,0}$ such that $\det\psi_{u,0}^+=1$ and $\psi_{u,0}^-=\psi_{u,0}^+\varepsilon$.
Since $\varepsilon^\sigma=\varepsilon$ and 
$\rho_{u,0}^{\sigma_3}=\rho_{3u,0}$, and together with
the extension with trivial determinant is unique, we obtain $(\psi_{u,0}^\pm)^{\sigma_3}=\psi_{3u,0}^\pm$.

Suppose that  $z\neq0$. If $z\equiv-z\pmod b$, then $z=0$ because $b$ is odd, this is a contradiction. Thus $\rho_{u,-z}\neq \rho_{u,z}$. However, every defect zero character of $N$ is writhed as $\psi_{u,z}=\Ind_U^{U\rtimes \langle \tau\rangle}(\rho_{u,z})$. So $\overline{\psi_{u,z}}=\Ind_S^{S\rtimes \langle \tau\rangle}(\overline{\rho_{u,z}})=\Ind_S^{S\rtimes \langle \tau\rangle}(\rho_{u,z}^\tau)=\psi_{u,z}$.  Since the induction is compatible with the action of $\sigma_3$,
we obtain $\psi_{u,z}^{\sigma_3}=\psi_{3u,3z}$.

For the $d_8$-family, we give the parametrization of $\dz_3(S)$ as in \eqref{3.21}. Let $c=(q^4-q^2+1)/3$. Similar to the $d_7$-family, starting from the \(d_8\)-family of ${\rm GU}_3(q^2)$, whose parametrization is $t\sim q^4t$ and  $t\in{\mathbb Z}/(q^6+1)$.
Restriction to \({\rm SU}_3(q^2)\), the parametrization can be taken as $k\sim q^4k$ and $k\in{\mathbb Z}/(q^4-q^2+1)$.
Here, \(k=0\) does not yield the required irreducible characters, and the two special parametrization $k=c$ and $k=2c$ give  degree-\(d_8/3\) characters. As for the $d_7$-family, the central element \(eI_3\) acts by the scalar \(e^k\), so the condition for descent to \(S\) remains
$3\mid k$. However, \(3\nmid c\), so neither of the two special 
parametrization \(c,2c\) satisfies this condition. 
We then write $k=3v$, where $v\in\mathbb Z/c$.
Excluding \(k=0\), we define
\begin{equation}\label{3.21}
 \{\rho_v|\ 0\ne v\in\mathbb Z/c,\ v\sim q^4v\}.   
\end{equation}
These characters have degree $\rho_v(1)=(q^2-1)(q^2+1)^2$ and 
are $3$-defect zero.

Suppose that $\rho_v^\tau=\rho_v$. Then 
$q^{4i}v\equiv-v\pmod c$ for some $i\in\{0,1,2\}$. Multiplication by $q^{4i}$ three times in succession yields
$q^{12i}v\equiv-v\pmod c$. But \(q^{12}\equiv1\pmod c\), so
$v\equiv-v\pmod c$. Since $c$ is odd, we obtain $v=0$, which is a contradiction. Thus $\rho_v\neq \rho_{-v}$.
Similar to the $d_7$-family, every defect zero character of $N$ 
is write as $\psi_{v}=\Ind_S^{S\rtimes \langle \tau\rangle}(\rho_{v})$. Then $\overline{\psi_v}=\Ind_S^{S\rtimes \langle \tau\rangle}(\overline{\rho_{v}})=\Ind_{S}^{S\rtimes  \langle \tau\rangle}(\rho_{v}^\tau)=\psi_v$. Similarly, for $\sigma_3$ as above, since $\rho_v^{\sigma_3}=\rho_{3v}$, we have $\psi_{v}^{\sigma_3}=\psi_{3v}$.

Finally, by \cite[Theorem 5.1]{Malle} and its proof, the action of field automorphisms on these defect zero characters is satisfied. We complete the proof.
\end{proof}

For all positive-defect $3$-blocks, applying \cite[Theorem 5.1]{Malle}, Lemmas \ref{eq:nine}-\ref{lem:3.20} and their proofs, we can match global characters and weight characters with the same parameterization, and this is compatible with the actions of Galois automorphisms and field automorphisms.
Together with the case of the defect zero $3$-block, we obtain
\begin{proposition}\label{pro:3.20}
There exists an $(\HH_3\times D)$-equivalent, preserving blocks, bijection between $\IBr_3(G)$ and $\WW_3(G)$.
\end{proposition}
\begin{proof}
For the principal $3$-block $B_0$, \cite[Theorem 5.1]{Malle} gives a bijection $\Omega_0$ between $\IBr_3(B_0)$ and the set of conjugacy classes of weights with radical group $R$ in Table \ref{tab:3}. In fact, for $a=1$ and $a\geq 2$, respectively,
the nine weights can be arbitrarily paired with the nine Brauer characters because they are $(\HH_3\rtimes D)$-invariant by Lemmas \ref{eq:nine} and \ref{lem:3.15}.

For a non-principal positive-defect $3$-blocks $B$, as in \eqref{Cs}, the centralizer \(\CC_{G}(s)\) is a torus 
or a direct product of \(\SL_2(q^2)\) and a tours.
For global part, we denote the corresponding one or two Brauer characters by $\eta_{s,1}$, or $\eta_{s,1}$ and $\eta_{s,\St}$, where  $\St$ is a unipotent character of \(\CC_{G}(s)\).
Moreover, by Lemma \ref{lem:3.19}, the action of $\HH_3\times D$ by acting on the semisimple $3'$-element $s$. 

For local part, according to different types of  $\CC_G(s)$, we have the following four cases.
If $\CC_G(s)\cong C_{q^2+1}\times\SL_2(q^2)$ or
$C_{q^2+1}\times C_{q^2+1}$, then by the proof of Lemma \ref{lem:3.20}, the weight characters are denoted as 
$\psi_{s}^+$ and $\psi_s^-$,  or $\psi_s$ (see \eqref{first}), respectively. In the former, we can set
$\eta_{s,1}\longmapsto\overline{(R,\psi_s^+)}$ and
$\eta_{s,\St}\longmapsto\overline{(R,\psi_s^-)}$.
In the latter, we set $\eta_{s,1}\mapsto \overline{(R,\psi_s)}$.

If $\CC_G(s)\cong C_{q^2-1}\times\SL_2(q^2)$, then when $a=1$, we can choose $\psi_{u,0}^+$ and $\psi_{u,0}^-$, which are two extensions of $\rho_{u,0}$(see \eqref{3.20}), as the weight characters. In the case, since $s$ corresponds to $u$, we set
$\eta_{s,1}:=\eta_{u,1}\longmapsto\overline{(R,\psi_{u,0}^+)}$ and
$\eta_{s,\St}:=\eta_{u,\St}\longmapsto\overline{(R,\psi_{u,0}^-)}$.
When $a\geq2$, we can choose $\psi_{\kappa,1}$ and $\psi_{\kappa,\varepsilon}$(see \eqref{3.19})  as the weight characters.
Note that \(\kappa:=\kappa_\lambda\) is the character of \(\SL_2(q^2)\) labelled by \(\lambda\sim -\lambda\) by the proof of $\eqref{3.19}$. Then $s$ corresponds to $\lambda$, and then we set $\eta_{s,1}:=\eta_{\lambda,1}\mapsto \overline{(R,\kappa_\lambda\otimes 1)}$ and
$\eta_{s,\St}:=\eta_{\lambda,\St}\mapsto \overline{(R,\lambda\otimes\varepsilon)}$.

If $\CC_G(s)\cong C_{q^2-1}\times C_{q^2+1}$, then when 
$a=1$, the weight character is chosen as 
$\psi_{u,z}=\Ind_S^{S\rtimes\langle\tau\rangle}(\rho_{u,z})$  which is the induction of $\rho_{u,0}$ (see \eqref{3.20} when $z\neq0$). In the case,  $(s,1)$ corresponds to \((u,z)\sim(-u,z)\sim(u,-z)\), and we set 
$\eta_{s,1}:=\eta_{u,z}\mapsto \overline{(R,\psi_{u,z})}$.
When $a\geq 2$, the weight character is chosen as 
$\psi_{\kappa_\lambda,\beta}$(see \eqref{3.19}) with $\beta=\rho_j\in\Irr(C_b\rtimes C_2)$. 
In the case,  $(s,1)$ corresponds to \((\lambda,j)\sim(-\lambda,j)\sim(\lambda,-j)\), and we set 
$\eta_{s,1}:= \eta_{\lambda,j}\mapsto \overline{(R,\psi_{\kappa_\lambda,\rho_j})}$.

If  $\CC_G(s)\cong C_{q^4-q^2+1}$, the weight character is chosen as  $\psi_{v}=\Ind_S^{S\rtimes\langle\tau\rangle}(\rho_{v})$, which is the induction of $\rho_{v}$ (see \eqref{3.21}) for $v\sim \pm q^{4i}v$ and $i=0,1,2$.
We set $\eta_{s,1}:=\eta_{v}\mapsto \overline{(R,\psi_{v}})$. This definition is independent of the choice of representative for \(v\) because $\psi_v=\psi_{\pm q^{4i}v}$.

For all non-principal positive-defect  $3$-blocks, by Lemma \ref{lem:3.20} and its proof, the action of
$\HH_3\times D$ on the parameterization of the weight characters
corresponds to its action on the element $s$. 
Moreover, the corresponding selection of the above local parameters does correspond the same block globally, which
is guaranteed  by \cite[Theorem 5.1]{Malle}.

Finally, by all correspondences as above, and together with the case of the defect zero blocks, we obtain the equivalent bijection. We complete the proof.
\end{proof}

Summarizing this Section, we have proven the following result.
\begin{theorem}\label{thm:2F_4}
Let $G={}^2F_4(2^{2m+1})$ with $m\geq 1$ and $p$ be odd. Then $G$ satisfies the inductive NAW condition at the prime $p$.
\end{theorem}
\begin{proof}
The cyclic Sylow case has been settled in Proposition \ref{lem:3.2}. When $p\geq5$ and Sylow $p$-subgroups are non-cyclic, by Theorem \ref{equ:bijection when p>=5}, we obtain the equivalent bijection. In the case, for any $\eta\in\IBr_p(G)$, we choose $A=G\rtimes D_\eta$ and $\widetilde A=G\rtimes D_{\eta^{\HH_p}}$. Note that $\eta$ is either real-valued, unipotent, or comes from one of the families $\chi_{23},\chi_{24},\chi_{43},\chi_{44},\chi_{47},\chi_{48}$.
So by Lemma \ref{lem:real}, Lemmas \ref{lem:unipotent}-\ref{lem: The non-split Suzuki-centralizer families}, $\eta$ has an equivalent
extension. By \cite[Lemma~10.7(ii)]{FFZ} and together with the map \eqref{eq:basic}, we obtain the invariant extensions of the global Brauer characters. All equivalent extensions of the local characters are obtained by Proposition \ref{equ:local extention p>=5}.  All hypotheses of Proposition \ref{prop:criterion} are satisfied. 

For $p=3$, by Proposition \ref{pro:3.20}, the equivalent bijection of Proposition \ref{prop:criterion} is satisfied, and  the resulting bijection preserves all blocks.  For any $\eta\in\IBr_3(G)$, we choose $A=G\rtimes D_\eta$ and $\widetilde A=G\rtimes D_{\eta^{\HH_3}}$. For the the principal $3$-block, by Lemmas \ref{eq:nine} and \ref{lem:3.15}, all global Brauer and local characters are $(\HH_3\rtimes D)$-invariant. Applying Lemma \ref{lem:h3}, we obtain the required extensions. For all non-principal positive-defect $3$-blocks, by Lemmas \ref{lem:3.19} and \ref{lem:3.20}, all global Brauer and local characters are real-valued. Applying  Lemma \ref{lem:real}, we obtain the required extensions. For the defect-zero $p$-block, this is obvious. 
All hypotheses of Proposition \ref{prop:criterion} are satisfied. We are done.
\end{proof}

\section{Remaining individual groups}\label{section 5}
In this Section, we will handle some small groups that were mentioned in Section \ref{section 2.3}. That is, we consider the derived subgroups ${}^2G_2(3)'$ and ${}^2F_4(2)'$. For the simple group ${}^2G_2(3)'$, which is isomorphic to  $\PSL_2(8)$. Then this case is handled by \cite[Proposition 10.8]{FFZ}.

From now on, we write \(G={}^2F_4(2)'\) and \(L={}^2F_4(2)\). Then $L\cong \Aut(G)$ and we identify \(G\) with \(\Inn(G)\). Note that the order $|G|$ has the  odd prime factors $p=3,5,13$ by \cite{T}. Moreover,  $\Out(G)=L/G\cong C_2$ which is only generated by an involution, and  $\CC_L(G)=1=\ZZ(L)$.
Since the Schur multiplier of $G$ is trivial, its own universal \(p'\)-covering group.  

We next consider $p=3$ and $p=5$, respectively.

\subsection{Verification at the prime  
\texorpdfstring{$p=3$}{p=3}}\label{sec:p3}
Let $P=P_3\in\Syl_3(G)$. Then by \cite{T}, $P\cong3_{+}^{1+2}$ is the  extraspecial group of order $27$, $\ZZ(P) = P'\cong  C_3$ and $\ZZ(P)/P\cong {\mathbb F}_2^3$. Moreover,  $\NN_G(P)\cong P\rtimes D_8$ has order $216$ and $\NN_L(P)\cong P\rtimes SD_{16}$ has order $432$.  The quotient $SD_{16}$ embeds in $\GL_2(3)$ via the action on $P/\ZZ(P)$. Let
\begin{equation}\label{p3:eq:SDmat}
u=\begin{pmatrix}1&1\\-1&1\end{pmatrix}\ \text{and}\
v=\begin{pmatrix}1&0\\0&-1\end{pmatrix}\ \text{over }\FF_3
\end{equation}
with \(u^8=v^2=1\) and \(vuv^{-1}=u^3\). 
Then \(D_8\cong\langle u^2,v\rangle\) and
\(SD_{16}\cong \langle u,v \rangle\).

\subsubsection{Global characters and two actions}
The 3-defect-zero ordinary characters of \(G\) are given by the ordinary characters table (see \cite{T}) and the \(3\)-decomposition matrix (see \cite[Block.1]{R}), we use their notation,
\begin{equation}\label{p3:eq:dzS}
 \dz_3(G)=\{\chi_4,\chi_5,\chi_9,\chi_{10},\chi_{11},
 \chi_{15},\chi_{16},\chi_{17},\chi_{20}\}.
\end{equation}
with degrees $27,27,351,351,351,675,702,702,1728$. 
The principal \(3\)-block $B_0$ of $G$ has nine Brauer characters, with labels and degrees(see \cite{R})
\begin{equation}\label{p3:eq:principal}
 \IBr_3(B_0)=\{\varphi_1,\varphi_2,\varphi_3,
 \varphi_6,\varphi_7,\varphi_8,\varphi_{12},
 \varphi_{16},\varphi_{17}\}
\end{equation}
with degrees $1,26,26,77,124,124,572,1099,1099$.
The other nine Brauer characters are exactly the restrictions \(\chi^\circ\) with \(\chi\in\dz_3(G)\) in \eqref{p3:eq:dzS}.   

Using the decomposition matrix in \cite{R} and the ordinary character values and class-power maps in \cite{T,OrdData}, one reconstructs all Brauer values on the eighteen \(3\)-regular classes. The third-power map $\pi_3$ is displayed in Appendix \ref{app:maps}. The fields generated by the values of the nine characters in $B_0$ are respectively $\QQ$ for $\varphi_1,\varphi_6,\varphi_{12}$, $\QQ(\sqrt{-2})$ for 
$\varphi_2,\varphi_3$, and $\QQ(\sqrt{13})$ for $\varphi_7,\varphi_8,\varphi_{16},\varphi_{17}$.

\begin{lemma}\label{equ}
All Brauer characters in the block $B_0$ are \(\mathcal H_3\)-invariant. Moreover, $\HH_3$ and $\Out(G)$ act respectively on $\IBr_3(G)$ by the non-trivial permutations
$\tau=(\varphi_4\ \varphi_5)(\varphi_{10}\ \varphi_{11})
 (\varphi_{14}\ \varphi_{15})$ and
$\alpha=(\varphi_2\ \varphi_3)
 (\varphi_7\ \varphi_8)
 (\varphi_{14}\ \varphi_{15})
 (\varphi_{16}\ \varphi_{17})$,
and other characters are $(\HH_3\times \Out(G))$-invariant.
\end{lemma}
\begin{proof}
For $\sigma_3\in\HH_3$, since $\sqrt{-2}=\zeta_8+\zeta_8^3$, we obtain $\sigma_3(\sqrt{-2})=\zeta_8^3+\zeta_8^9=\sqrt{-2}$.
Let $\zeta_{13}$ is $13$-th root of unity.
Define $\chi(a)=\left(\frac{a}{13}\right)$,
where $\left(\frac{a}{13}\right)$ is the Legendre symbol. The corresponding quadratic Gauss sum is
$\tau_{13}=\sum_{a=1}^{12}\chi(a)\zeta_{13}^{a}$.
Since $13\equiv 1\pmod{4}$, the classical formula for quadratic Gauss sums gives $\tau_{13}=\sqrt{13}$.
Thus, $\tau_{13}$ generates 
$\mathbb{Q}(\tau_{13})=\mathbb{Q}(\sqrt{13})$.
Note that \(3\equiv4^2\pmod{13}\), so 
\(\sigma_3(\tau_{13})=\tau_{13}\). Thus we have
\(\sigma_3(\sqrt{13})=\sqrt{13}\).    
By \cite[Tables]{T}, the action of $\sigma_3$ on the set $\IBr_3(G)$ by the non-trivial permutation is $\tau=(\varphi_4\ \varphi_5) (\varphi_{10}\ \varphi_{11})(\varphi_{14}\ \varphi_{15})$, and  all other \(\varphi_i\) are $\HH_3$-fixed. Moreover, these three pairs occur in the $3$-defect-zero part with degrees \(27,351,702\), respectively.

Let \(1\neq \alpha\in \Out(G)\cong C_2\). By the stored class fusion of $G$ into $L$, which refers to a mapping from the conjugate class set of $G$ to the conjugate class set of $L$ and is recorded in \cite{OrdData} or Appendix \ref{app:maps}, and  together with two decomposition matrices of $G$ and $L$(see \cite{T,R}), we have the permutation
$\alpha=(\varphi_2\ \varphi_3) (\varphi_7\ \varphi_8)(\varphi_{14}\ \varphi_{15})(\varphi_{16}\ \varphi_{17})$. 
Thus in the block $B_0$ has three  $\alpha$-fixed points $\varphi_1,\varphi_6,\varphi_{12}$ and three $\alpha$-orbits
$\{\varphi_2\ \varphi_3\}, \{\varphi_7\ \varphi_8\}, \{\varphi_{16}\ \varphi_{17}\}$.
We complete the proof.
\end{proof}

\subsubsection{Local characters and two actions} 
By \cite[Theorem 5.1 and its proof]{Malle} and \eqref{p3:eq:dzS}, there are nine  $3$-weights with trivial radical group. These nine are separate $3$-defect-zero blocks.
Moreover, there are nine $3$-weights with non-trivial radical group $R$ as follows. These results can be founded in \cite{T}, for the completeness of the paper, we still provide a proof.
\begin{lemma}\label{p3:prop:radical}
There are exactly three non-trivial \(G\)-conjugacy classes of \(3\)-radical subgroups, represented by $R_+, R_-, P$, where \(R_\pm\cong C_3\times C_3\) satisfies
 $\NN_G(R_\pm)/R_\pm\cong\GL_2(3)$, and the non-trivial outer automorphism exchanges the two \(G\)-classes \([R_+]\) and \([R_-]\).
\end{lemma}
\begin{proof}
Since \(\NN_G(P)/P\cong D_8\) is a \(3'\)-group, \(\OO_3(\NN_G(P))=P\), and then $P$ is radical.
By \cite{A}, the ATLAS has a single class \(3X\) of elements of order \(3\), with \(|\CC_G(z)|=108\) for \(z\in3X\).
Since \(Z(P)\) is characteristic in \(P\),
\(\NN_G(P)\leq\NN_G(Z(P))\).  On the other hand,
\(\NN_G(\langle z\rangle)/\CC_G(z)\leq C_2\). Thus
\[
 216=|\NN_G(P)|\leq|\NN_G(Z(P))|\leq2\times108=216.
\]
Thus \(\NN_G(Z(P))=\NN_G(P)\), whose normal \(3\)-core is \(P>Z(P)\). All subgroups of order \(3\) are conjugate, so they are not radical.

The GAP maximal-subgroup list (see \cite{T}) contains \(L_3(3):2\leq G\). Inside its normal subgroup \(L_3(3)=\SL_3(3)\). Let \(Q=R_0:\GL_2(3)\) with \(R_0\cong C_3\times C_3\), and then \(|Q|=432\).  Take its Sylow \(3\)-subgroup \(P\) so that \(R_0\leq P\).  Since \(R\leq\NN_G(R_0)\) and  $|\NN_G(R_0)|\leq 432$, this forces \(\NN_G(R_0)=Q\).  Thus \(R_0\) is radical because
\(\OO_3(\GL_2(3))=1\).

Note that the action of \(\NN_L(P)/P\cong SD_{16}\) on the four order-\(9\) subgroups of \(P\) is transitive by \eqref{p3:eq:SDmat}.  Thus they
are radical and have normalizer quotient \(\GL_2(3)\).
Note that \(\NN_G(P)/P=D_8\) has exactly two orbits on
the four subgroups.  These remain distinct under \(G\)-conjugation: if \(R^g=R'\) with \(R,R'\leq P\), then both \(P^g\) and \(P\) are Sylow \(3\)-subgroups of \(\NN_G(R')\).  Sylow conjugacy inside \(\NN_G(R')\) supplies \(h\in\NN_G(R')\) with \(P^{gh}=P\) and \(R^{gh}=R'\).  Therefore \(gh\in\NN_G(P)\), so the two subgroups belonged to the same \(D_8\)-orbit already. The transitivity of \(\NN_L(P)\) shows that the outer involution interchanges these two \(G\)-classes.
\end{proof}

\begin{lemma}\label{equ1}
Let $R$ be a radical subgroup in Lemma \ref{p3:prop:radical}.
Then all defect zero characters of $\NN_G(R)/R$  are rational-valued.
\end{lemma}
\begin{proof}
If \(R=R_\pm\), then \(\NN_G(R)/R=\GL_2(3)\) by Lemma \ref{p3:prop:radical}, whose irreducible ordinary degrees are $1,1,2,2,2,3,3,4$ by \cite{T}.
Since \(|\GL_2(3)|_3=3\),  its two characters of degree
\(3\) are \(3\)-defect zero. They are inflated from the
two degree-\(3\) characters of \(\operatorname{PGL}_2(3)\cong S_4\). Thus they are rational-valued. If \(R=P\), then \(\NN_G(R)/R=D_8\), which is a \(3'\)-group and has five ordinary irreducible characters with degrees \(1,1,1,1,2\) by \cite{T}. Moreover, they are $3$-defect zero characters and rational-valued. 
\end{proof}

Let \(1\neq \alpha\in \Out(G)\). For $R=P$, write
\(D_8=\langle r=u^2,s=v\rangle\) by using \eqref{p3:eq:SDmat}. Denote its linear characters by
\(\lambda_{\epsilon,\delta}(r)=\epsilon\),
\(\lambda_{\epsilon,\delta}(s)=\delta\) for
\(\epsilon,\delta\in\{\pm1\}\), and its degree-\(2\) character by \(\rho\).  
With the convention \(d^u=u^{-1}du\), conjugation sends
\(r\mapsto r\), \(s\mapsto rs\).  Thus, $\alpha$ fixes
 $\lambda_{+,+},\ \lambda_{+,-},\ \rho$ and interchanges \(\lambda_{-,+}\) with \(\lambda_{-,-}\).
For \(R=R_\pm\), we denote the two degree-\(3\) characters of \(\GL_2(3)\) by \(\tau_{\pm,0}\) and \(\tau_{\pm,1}\) such that $\alpha$ sends \(\tau_{+,j}\) to \(\tau_{-,j}\), \(j=0,1\). 
Also, by Lemma \ref{p3:prop:radical}, $\alpha$ exchanges 
\(Q_+\)-weights and \(Q_-\)-weights. Consequently, the nine non-trivial weights have the same $\alpha$-orbit structure as \(\IBr_3(B_0)\).
\begin{table}[ht]
\centering
\caption{One correspondence on $B_0$}\label{tab:6}
\begin{tabular}{@{}ll@{}}
\toprule
$3$-Brauer characters&$3$-Weights\\
\midrule
\(\varphi_1\)&\(\overline{(P,\lambda_{+,+})}\)\\
\(\varphi_6\)&\(\overline{(P,\lambda_{+,-})}\)\\
\(\varphi_{12}\)&\(\overline{(P,\rho)}\)\\
\(\varphi_2,\ \varphi_3\)&\(\overline{(P,\lambda_{-,+})},\ \overline{(P,\lambda_{-,-})}\)\\
\(\varphi_7,\ \varphi_8\)&\(\overline{(Q_+,\tau_{+,0})},\ \overline{(Q_-,\tau_{-,0})}\)\\
\(\varphi_{16},\ \varphi_{17}\)&\(\overline{(Q_+,\tau_{+,1})},\ \overline{(Q_-,\tau_{-,1})}\)\\
\bottomrule
\end{tabular}
\end{table}

\begin{proposition}\label{p3:prop:equivariant}
There exists an 
\((\HH_3\times\Aut(G))\)-equivalent bijection 
\(\Omega:\IBr_3(G)\longrightarrow\mathcal W_3(G)\).
\end{proposition}
\begin{proof}
Define \(\Omega\) on the defect-zero blocks by
$\chi^\circ\mapsto\overline{(1,\chi)}$ for every $\chi\in\dz_3(G)$ in \eqref{p3:eq:dzS}. The map commutes with both Galois automorphisms and outer automorphisms. On the block \(B_0\), the action of $\HH_3$ is trivial on both sides by Lemmas \ref{equ} and \ref{equ1}. By two permutations
in Lemma \ref{equ}, and together with actions on the weight characters, we  establish the required bijection as in Table \ref{tab:6}. 
\end{proof}

\subsubsection{Compatible extensions of global and local characters}
For the extensions of Brauer characters and weight characters in the above bijection. Note that these characters are $\HH_3$-invariant, by Lemma \ref{lem:h3},
they have the required extensions, respectively.
However, we provide the specific character extensions.
\begin{proposition}\label{thm:p3}
Let $\eta\in\IBr_3(G)$ and $\Omega(\eta)=\overline{(R,\psi)}\in\WW_3(G)$.
Then $\eta$ and $\psi^\circ$ satisfy Proposition \ref{prop:criterion}.
\end{proposition}
\begin{proof}
Choose \(\AT=L_{\eta^{\HH_3}}\) and \(A=L_\eta\). Note that $\eta$ comes from restriction of \eqref{p3:eq:dzS} and \eqref{p3:eq:principal}, and that $\psi$ comes from Table \ref{tab:6}. First, for the fixed three characters \(\varphi_1,\varphi_6,\varphi_{12}\) in the principal block $B_0$, their degrees are \(1,77,572\) by \eqref{p3:eq:principal}. Let \(\widehat\chi_6\) and \(\widehat\chi_{14}\) be ordinary rows \(8\) and \(21\) of \(L\) in \cite{OrdData}. They have degrees \(78\) and \(650\), and their restriction  to $G$ are \(\chi_6\) and \(\chi_{14}\). Moreover,  \(\widehat\chi_6\) and \(\widehat\chi_{14}\) are rational-valued. The subscripts here are their \(G\)-labels. The pertinent rows of \cite[Block. 1]{T} are
\[
 \chi_{6,0}^{\circ}=\varphi_{1,0}+\varphi_{6,0},\quad
\chi_{14,0}^{\circ}=\varphi_{1,0}+\varphi_{6,0}+\varphi_{12,0}.
\]
Here, the double subscripts are the labels for \(L\).
Choose the signs of the extensions compatibly with these rows. Then 
$\widehat\varphi_1=1_L,\ \widehat\varphi_6=\widehat\chi_6^{\circ}-1_L,\ 
\widehat\varphi_{12}=\widehat\chi_{14}^{\circ}-\widehat\chi_6^{\circ}$
are irreducible Brauer extensions of the required degrees. In particular, irreducibility of these
differences uses the decomposition matrix of $L$, not their degrees alone. All three are rational-valued, and thus
\(\HH_3\)-invariant.  

For the local side, in the case, 
\(\NN_L(P)/P\cong SD_{16}=\langle u,v\rangle\) by \eqref{p3:eq:SDmat}. For two fixed linear characters  
$\lambda_{+,+}$ and $\lambda_{+,-}$ of $D_8$, they extend to $SD_{16}$ by
 $\widehat\lambda_{+,+}(u)=1,\ \widehat\lambda_{+,+}(v)=1,\ 
 \widehat\lambda_{+,-}(u)=1,\ \widehat\lambda_{+,-}(v)=-1$.
These extensions are rational-valued. For the degree-$2$ character \(\rho\), we take \(\kappa\in\Irr(\langle u\rangle)\) with \(\kappa(u)=\zeta_8\), and let
 $\widehat\rho=\Ind_{\langle u\rangle}^{SD_{16}}(\kappa)$. 
The stabilizer of \(\kappa\) is \(\langle u\rangle\), so
\(\widehat\rho\) is irreducible of degree $2$. Its restriction to \(D_8=\langle u^2,v\rangle\) is \(\rho\). Since \(\kappa^{\sigma_3}=\kappa^3=\kappa^v\), induction of $\kappa$ 
gives \(\widehat\rho^{\sigma_3}=\widehat\rho\).
Thus, the three local extensions, inflated to \(\NN_L(P)\), are \(\HH_3\)-invariant.

For these pairs in the former three rows of Table \ref{tab:6}, \(\AT=A=L\). Conjugation by \(\NN_L(P)\) fixes the global and local extensions. Together with Galois invariance this makes them stable under the full combined stabilizer, as required.

For every member of 
\( \{\varphi_2,\varphi_3,\varphi_7,\varphi_8,
\varphi_{16},\varphi_{17}\}\) in the latter three rows of Table \ref{tab:6}, it's \(\HH_3\)-invariant and has an
outer automorphism orbit of length two by Lemma \ref{equ}. Thus \(\AT=A=G\). The corresponding weight has the same stabilizer by Proposition~\ref{p3:prop:equivariant}.
Take \(\widehat\eta=\eta\) and \(\widehat\psi=\psi^\circ\), as required. Finally, for the nine defect-zero characters, it is  obvious. We complete the proof.
\end{proof}

\subsection{Verification at the prime  
\texorpdfstring{$p=5$}{p=5}}\label{sec:p5}
Let \(P:=P_5\in\Syl_5(G)\). Then $P\cong C_5\times C_5$, $\NN_G(P)\cong P\rtimes 4A_4$ and $\NN_L(P)\cong P\rtimes 4S_4$(see \cite{T}).
Write \(E=\NN_G(P)/P\) and \(\et=\NN_L(P)/P\). 
Then $E\cong 4.A_4$ and $\et\cong 4.S_4$.  In the standard affine realization, we obtain
\begin{equation}\label{p5:eq:central-product}
 E\cong C_4\circ\SL_2(3),\
 \ZZ(E)=\ZZ(\et)\cong C_4,\
 E/\ZZ(E)\cong A_4,\
 \et/\ZZ(\et)\cong S_4.
\end{equation}
The central product identifies the two central involutions.
The table-library entries in \cite{T} identify the affine normalizer $\NN_G(P)$  with \texttt{PrimitiveGroup(25,18)}, and $\NN_L(P)$  with \texttt{PrimitiveGroup(25,19)}, respectively.
Record the local inclusion through a  class fusion
in \cite{T}.

\subsubsection{Global characters and two actions}
All ordinary characters \(\chi_i\) and Brauer characters \(\varphi_j\) below are those of \cite{T}. 
The decomposition matrix of the principal block is reproduced in \cite[Block.1]{T}. Values of ordinary characters and the class fusion map \(G\hookrightarrow L\) are taken from the
named table records in \cite{OrdData} and Appendix \ref{app:maps}. 

There are one positive-defect \(5\)-block, the principal block \(B_0\), and six defect-zero blocks.  Moreover, the latter have
\begin{equation}\label{p5:dz}
 \dz_5(G)=\{\chi_7,\chi_8,\chi_{14},\chi_{15},\chi_{18},\chi_{19}\}
\end{equation}
with degrees $300,325,650,675,1300,1300$.
Their restrictions to $5$-regular elements are
$\varphi_9$, $\varphi_{10}$, $\varphi_{17}$, $\varphi_{18}$, $\varphi_{19}$, $\varphi_{20}$.
Here, the labels and degrees are directly from the top table of \cite{T}. Moreover, by \cite[Block.1]{T},  $\IBr(B_0)$ has fourteen irreducible Brauer characters
\begin{equation*}
\IBr_5(B_0)=\{\varphi_1,\varphi_2,\varphi_3,\varphi_4,\varphi_5,\varphi_6,
\varphi_7,\varphi_8,\varphi_{11},\varphi_{12},\varphi_{13},\varphi_{14},\varphi_{15},\varphi_{16}\}
\end{equation*}
with degrees $1,26,26,27,27,78,109,109,351,351,460,460,593,593$.
\begin{lemma}\label{lem:p=5 global characters}
All Brauer characters in defect-zero blocks are $(\HH_5\times \Out(G))$-invariant. Moreover, $\HH_5$ and $\Out(G)$ act respectively on the set $\IBr_5(B_0)$ by the non-trivial permutations 
$\tau=(\varphi_2\ \varphi_3)(\varphi_7\ \varphi_8)
 (\varphi_{13}\ \varphi_{14})(\varphi_{15}\ \varphi_{16})$
and
$\alpha=(\varphi_2\ \varphi_3)(\varphi_7\ \varphi_8)
 (\varphi_{13}\ \varphi_{14})(\varphi_{15}\ \varphi_{16})
 (\varphi_{19}\ \varphi_{20})$,
and other characters are $(\HH_5\times\Out(G))$-invariant.
\end{lemma}
\begin{proof}
Let \(D\) be the  decomposition matrix of $B_0$ as in \cite[Block. 1]{T}. Then ${\rm rank}(D)=14$.
Restrict the sixteen ordinary rows in $B_0$ to the twenty \(5\)-regular classes, obtaining a matrix \(Y\). Write \(X\) for the table of Brauer characters for the principal block $B_0$. Then  $DX=Y$ and such 
$X$ is unique. The fifth-power map on classes computes $\sigma_5$ on the rows of $X$(see Appendix \ref{app:maps}).

On the other hand, we consider the class fusion map 
$G\hookrightarrow L$. Since $L/G\cong \Out(G)$,
two \(G\)-conjugacy classes \(c_1,c_2\) are fused in \(L\) if and only if $c_1=\alpha(c_2)$ for some $\alpha\in L/G$.
Thus this compute the action of the outer automorphism $\alpha\neq 1$. The input maps are written explicitly in  Appendix \ref{app:maps}. From \cite{T,OrdData}, the calculation gives  $\tau
=(\varphi_2\ \varphi_3)(\varphi_7\ \varphi_8)(\varphi_{13}\ \varphi_{14})(\varphi_{15}\ \varphi_{16})$ and $\alpha=(\varphi_2\ \varphi_3)(\varphi_7\ \varphi_8)(\varphi_{13}\ \varphi_{14})(\varphi_{15}\ \varphi_{16})(\varphi_{19}\ \varphi_{20})$, and other characters are fixed. 
Thus the two actions agree on $\IBr(B_0)$, with fixed set
$\{\varphi_1, \varphi_4,\varphi_5,\varphi_{11},\varphi_{12} \}$. They fix all six defect-zero characters under $\tau$, whereas $\alpha$ exchanges the two of degree $1300$.
We complete the proof.
\end{proof}

\subsubsection{Local characters and two actions}
When $p=5$, all \(5\)-radical subgroups of $G$ and some local group structures can be founded in \cite{T}. For the completeness of the paper, we still provide the corresponding  proofs.
\begin{lemma}\label{p5:prop:radical}
There is one non-trivial \(G\)-conjugacy classes of \(5\)-radical subgroups, represented by $P$.
\end{lemma}
\begin{proof}
Note that \(\OO_5(\NN_G(P))=P\) because \(P\) is a normal Sylow subgroup of its normalizer. 
Let \(Q=\langle z\rangle\leq P\) have
order \(5\). Since \(P\) is abelian, $P\leq\CC_G(Q)=\CC_G(z)$ which is of order $50$ because the ATLAS conjugacy-class table in \cite{A} has one class of elements of order $5$, labelled $5A$, with centralizer order $50$. Also, a group of order $50$ has a unique Sylow $5$-subgroup, so
\(\OO_5(\CC_G(Q))=P\). The centralizer is normal in \(\NN_G(Q)\), and then its \(5\)-core is characteristic. 
Thus $\OO_5(\NN_G(Q))=P>Q$. So \(Q\) is not radical. 
We are done.
\end{proof}

Therefore, every weight character comes either from $\dz_5(G)$ with trivial radical $5$-subgroup, or from
$\Irr(E)$ with non-trivial radical $5$-subgroup $P$. Since $E$ is a $5'$-group, all its irreducible characters are $5$-defect zero. For describe the weight characters of $G$, we need some group structures as follows. In \(\GL_2(5)\), we take two elements
\[
a=\begin{pmatrix}2&0\\0&3\end{pmatrix},\quad b=\begin{pmatrix}0&1\\4&0\end{pmatrix},
\]
and write $K=\langle a,b\rangle\cong Q_8$.
Then the affine realizations above can be taken as
\begin{equation*}
 \et=\NN_{\GL_2(5)}(K),\quad
 E=\{g\in\et|\det(g)\in\{1,4\}\}.
\end{equation*}
This is the pair of local actions on \(P\cong\FF_5^2\) used below. For example, let
\[
t=\begin{pmatrix}3&3\\4&1\end{pmatrix},\quad
c=2I_2,\quad
u=\begin{pmatrix}0&4\\3&0\end{pmatrix}.
\]
Then \(t\) has order \(3\), \(c\) has order \(4\), and
$\langle a,b,t\rangle\cong\SL_2(3)$. Thus
$E=\langle a,b,t,c\rangle$ and  $\et=\langle E,u\rangle$.
Here \(c\) is central and \(u\notin E\) because \(\det(u)=3\). The projective images of \(E\) and \(\et\) are \(A_4\) and \(S_4\). These matrix calculations establish the descriptions in \eqref{p5:eq:central-product}.

Since \(5\nmid|\et|\), ordinary irreducible characters and Brauer characters of \(E\) and \(\et\) are naturally identified. 
Let $V$ be the natural two-dimensional \(\mathbb F_5E\)-module. Its second symmetric $\Sym^2 V$ is a three-dimensional \(\FF_5E\)-module. Write \(\nu\) for the Brauer character of $V$, regarded also as the corresponding ordinary character. 
Let $\delta=\det|_E$ and $\theta=\operatorname{Br}(\Sym^2 V)$.
Here, $\delta$ is viewed as a Brauer linear character by lifting its eigenvalue to a root of unity, and it has order \(2\).  $\theta$ denotes the Brauer character corresponding to the \(\FF_5E\)-module $\Sym^2 V$. Finally, choose a linear character \(\lambda\) of order \(3\) with the kernel $\langle K,c\rangle$.
\begin{lemma}\label{p5:lem:localchars}
The following list consists of all irreducible characters of \(E\)
\begin{equation*}
 \begin{split}
 \Irr(E)={}&\{\delta^e\lambda^j|\ e\in\{0,1\},\ j\in\{0,1,2\}\}\\
 &\sqcup\{\nu\delta^e\lambda^j| \ e\in\{0,1\},\ j\in\{0,1,2\}\}\\
 &\sqcup\{\theta,\theta\delta\},
 \end{split}
\end{equation*}
and their degrees are \(1^6,2^6,3^2\), where the superscripts denote multiplicities.
\end{lemma}
\begin{proof}
We claim that \(V\) is absolutely  irreducible on \(K\). Take the standard basis 
$e_1,e_2$ of $V$. Then $ae_1=2e_1$, $ae_2=3e_2$, $be_1=4e_2$ and  $be_3=e_1$.
Suppose there exists a one-dimensional \(K\)-invariant subspace \(W\). Since \(W\) is \(a\)-invariant, and the two eigenvalues \(2,3\) of \(a\) are distinct, \(W\) must be one of  $\langle e_1\rangle$ and $\langle e_2\rangle$. However,
$b\langle e_1\rangle=\langle e_2\rangle$ and 
$b\langle e_2\rangle=\langle e_1\rangle$. This shows that two one-dimensional subspace is not \(K\)-invariant. This is a contradiction. Thus 
$V$ is  irreducible on \(K\). Note that this argument still holds over the algebraic closure, so $V$ remains absolutely irreducible, as claimed.
Since \(K\leq E\), $\nu\in \Irr(E)$ with $\nu(1)=2$,
and thus $\nu\delta^e\lambda^j\in\Irr(E)$.

Note that $\dim \Sym^2V=3$, and then take the basis $e_1^2, e_1e_2, e_2^2$. By the action of $a$ and $b$, we have three common eigenvectors
$v_1=e_1^2+e_2^2$, $v_2=e_1^2-e_2^2$, $v_3=e_1e_2$.
Thus ${\rm Res}_K^E(\Sym^2V)=\bigoplus_{i=1}^3 \langle v_i \rangle$. Since $K/\ZZ(K)\cong C_2\times C_2$, the restriction ${\rm Res}_K^E(\theta)$ is the sum of three non-trivial linear characters $\alpha_i$.  Since $tat^{-1}=b^{-1}$ and $tbt^{-1}=-ab$, \(t\)-conjugation cycles these three $\alpha_i$,
Suppose that  $0\neq W\subseteq \Sym^2V$ is an 
$E$-invariant subspace. Since $K\unlhd E$ and $5\nmid|K|$, $W$ is reducible on $K$.
Note that $W$ must be the direct sum of several of them. Due to the cyclic permutation of these three constituents, there can only be $W=\Sym^2V$.
Thus $\theta\in \Irr(E)$ with $\theta(1)=3$.
The same is true for $\theta\delta$.

Since $c=2I_2$, \(c\) is a scalar matrix on \(V \).
Under the lift of Brauer characters,  
$(\nu\delta^e\lambda^j)(c)=(-1)^e\nu(c)$ because $\delta(c)=-1$, $\lambda(c)=1$ and $\nu(c)\neq0$. On the other hand, the trace of $t$ is $-1$
and  its determinant is $1$, so the characteristic polynomial of $t$ is $x^2+x+1$. Its two eigenvalues, when lifted to the complex numbers, are $\omega,\omega^2$. Thus $\nu(t)=\omega+\omega^2=-1$.
Also, since $\delta(t)=1$ and $\lambda(t)=\omega$,
we obtain $\nu\delta^e\lambda^j)(t) =-\omega^j$.
Therefore, the six degree-$2$ characters $\nu\delta^e\lambda^j$ are pairwise distinct.
The six linear characters $\delta^e\lambda^j$ are distinct since \(\delta\) and \(\lambda\) have coprime orders. 

Finally, 
note that $6\cdot1^2+6\cdot2^2+2\cdot3^2=48=|E|$, so there cannot be any other irreducible characters.
\end{proof}

We now consider the actions of outer automorphisms and Galois automorphisms on the set $\Irr(E)$.
\begin{lemma}\label{lem:p=5 local characters}
The Brauer restrictions of these six characters $1,\delta,\nu,\nu\delta,\theta,\theta\delta$ are $(\HH_5\times(\et\setminus E))$-invariant. 
Moreover, the action of $\HH_5\times(\et\setminus E)$ on other characters induced by the map $\lambda\mapsto\lambda^{-1}$.
\end{lemma}
\begin{proof}
Note that the six  characters $1,\delta,\nu,\nu\delta,\theta,\theta\delta$ in Lemma \ref{p5:lem:localchars} are afforded by \(\FF_5\)-modules  $1,\delta, V, V\otimes \delta, \Sym^2 V, \Sym^2 V\otimes \delta$, respectively. Then their Brauer values are fixed by 
\(\sigma_5\in\HH_5\).

Note that conjugation by an element of \(\et\setminus E\) fixes the six characters. On the other hand, $\lambda^{\sigma_5}=\lambda^5=\lambda^{-1}$ because 
$\lambda$ has order $3$. An odd permutation in \(S_4\) inverts \(A_4/V_4\cong C_3\), so $\lambda^\alpha=\lambda^{-1}$ for any $\alpha\in \et\setminus E$. Also, $\delta$ and $\nu$ have order $2$. Therefore, we obtain that
 $(\delta^e\lambda^j)^{\sigma_5}=(\delta^e\lambda^j)^\alpha=\delta^e\lambda^{-j}$ and
 $(\nu\delta^e\lambda^j)^{\sigma_5}=(\nu\delta^e\lambda^j)^\alpha=\nu\delta^e\lambda^{-j}$.
We complete the proof.
\end{proof}

Thus, for all members of $\Irr(E)$, both actions have
six  $\sigma_5\alpha$-fixed points and four 
$\sigma_5\alpha$-orbits. Together with Lemma \ref{lem:p=5 global characters}, we have the following correspondence.
\begin{table}[htbp]
\centering
\caption{The correspondence between Brauer characters and weight characters}\label{tab:7}
\begin{tabular}{lll}
\toprule
Brauer characters & Degrees & Weight characters\\
\midrule
\(\varphi_1\) & 1 & \(1\)\\
\(\varphi_6\) & 78 & \(\delta\)\\
\(\varphi_4\) & 27 & \(\nu\)\\
\(\varphi_5\) & 27 & \(\nu\delta\)\\
\(\varphi_{11}\) & 351 & \(\theta\)\\
\(\varphi_{12}\) & 351 & \(\theta\delta\)\\
\((\varphi_2,\varphi_3)\) & \((26,26)\) & \((\lambda,\lambda^2)\)\\
\((\varphi_7,\varphi_8)\) & \((109,109)\) & \((\delta\lambda,\delta\lambda^2)\)\\
\((\varphi_{13},\varphi_{14})\) & \((460,460)\) & \((\nu\lambda,\nu\lambda^2)\)\\
\((\varphi_{15},\varphi_{16})\) & \((593,593)\) & \((\nu\delta\lambda,\nu\delta\lambda^2)\)\\
\bottomrule
\end{tabular}
\end{table}

\begin{proposition}\label{p5:prop:equivariant}
There exists an 
\((\HH_5\times\Aut(G))\)-equivalent bijection
\(\IBr_5(G)\longrightarrow\mathcal W_5(G)\).
\end{proposition}
\begin{proof}
Define \(\Omega\) on the defect-zero blocks by
$\chi^\circ\mapsto\overline{(1,\chi)}$ for every $\chi\in\dz_5(G)$ in \eqref{p5:dz}. For the principal block $B_0$, we choose the correspondence in Table \ref{tab:7}. We compare both actions in Lemma \ref{lem:p=5 global characters} with in Lemma \ref{lem:p=5 local characters},  both sides fix the first--six rows and interchange the entries in each of the four paired rows. The inner automorphisms act trivially on both sets of conjugacy classes. We are done.
\end{proof}

\subsubsection{Compatible extensions of global and local characters} 
The following verifications concern  extensions of Brauer characters and weight characters for the above bijection. 
\begin{proposition}\label{thm:p5}
Let $\eta\in\IBr_5(G)$ and $\Omega(\eta)=\overline{(R,\psi)}\in\WW_5(G)$.
Then $\eta$ and $\psi^\circ$ satisfy Proposition 
\ref{prop:criterion}.
\end{proposition}
\begin{proof}
Choose \(\AT=L_{\eta^{\HH_5}}\) and \(A=L_\eta\).
For the Sylow $5$-subgroup \(P\), the Frattini argument gives \(L=G\NN_L(P)\), so $L/G\cong \NN_L(P)/\NN_G(P)$, as required. 
Note that $\eta$ comes from restriction of $\dz_5(G)$ and $\IBr_5(B_0)$, and that $\psi$ comes from Lemma \ref{p5:lem:localchars}. 

For extensions of the fixed characters
$\varphi_1,\varphi_4,\varphi_5,\varphi_6,
\varphi_{11},\varphi_{12}\in\IBr_5(B_0)$. 
As in \cite[Block.1]{T}, 
$\varphi_j=\chi_j^\circ\ (j=1,4,5,6)$, 
$\varphi_{11}=\chi_{10}^\circ$ and $\varphi_{12}=\chi_{11}^\circ$. For these six ordinary characters $\chi_j\ (j=1,4,5,6)$, $\chi_{10}$, $\chi_{11}$, there are extensions to \(L\) with values in \(\QQ(i)\), 
where $i$ is the imaginary unit.  More precisely, 
the following rows of $L$  in the record \texttt{MOT("2F4(2)'.2",...)} of \cite{OrdData} are suitable as follows. These characters $\varphi_1,\varphi_4, \varphi_5, \varphi_6, \varphi_{11}, \varphi_{12}$ correspond 
to rows $1,4,6,8,16,18$ in $L$. Here, these row numbers in $L$ refer to the ordinary character table of \(L\) record. Their restrictions give characters of \(G\). 

Let \(\widehat\chi\) be any of extensions of the six chosen ordinary characters, and \(\eta={\rm Res}_G^L(\widehat\chi^\circ)\). 
By an ordinary representation of \(L\) and reduction modulo \(5\), this gives a $\FF_5$-representation of \(L\) whose restriction to \(G\) has the character \(\eta\), which is irreducible. Note that the above  Brauer restriction
is a simple module, so the representation of \(L\) is itself simple. Thus \(\widehat\chi^\circ\) is an extension of \(\eta\).
By the decomposition matrix of \(L\)(see \cite[Block~1, pp.1--2]{T}),  independently confirm these six characters.
Since \(i^5=i\), every such extension is \(\HH_5\)-invariant. In particular, all $\varphi_1,\varphi_4,\varphi_5,\varphi_6,\varphi_{11},\varphi_{12}$ have an \(\HH_5\)-invariant extension to \(L\), denoted by
$\widehat{\varphi_1}, \widehat{\varphi_4}, \widehat{\varphi_5},
\widehat{\varphi_6}, \widehat{\varphi_{11}}, \widehat{\varphi_{12}}$.

For the fixed characters $1,\delta,\nu,\nu\delta,\theta,\theta\delta\in\Irr(E)$, we will construct their extensions. It is easy to check that each also extends to \(\et\) by the corresponding module extensions. Using its natural module \(\widehat V=\FF_5^2\), its character \(\widehat\delta=\det|_{\widetilde E}\), and the six extended modules 
\begin{equation}\label{p5:eq:local extensions}
 1,\ \widehat\delta,\ \widehat V,\ \widehat V\otimes\widehat\delta,\
 \Sym^2\widehat V,\ \Sym^2\widehat V\otimes\widehat\delta.
\end{equation}
Since the six extended modules are defined over \(\FF_5\), they are also \(\HH_5\)-invariant. We choose \(\AT=A=L\). Match every global extension $\widehat\varphi_i$ with  local extensions in \eqref{p5:eq:local extensions}, inflated from \(\et\) to \(\NN_L(P)\). Therefore, both extensions are \((\HH_5\times \NN_L(P))\)-invariant by Lemmas \ref{lem:p=5 global characters} and \ref{lem:p=5 local characters}.

Let \(\eta\) be one of the last four rows in Table \ref{tab:7}. Then by Lemma \ref{lem:p=5 global characters}, $\eta^\alpha=\eta^{\sigma_5}\ne\eta$. It is easy to check that $L_\eta=G$ and $L_{\eta^{\HH_5}}=L$. Thus \(\AT=L\), but \(A=G\). In particular, an element combining an outer automorphism with an odd Frobenius power can fix $\eta$.
On the other hand, note that $A/G=1=\NN_A(P)/\NN_G(P)$, so
there is no non-trivial extension, as required.
Finally, for the six defect-zero characters, $L_{\chi_j}=L$ for $j=7,8,14,15$, and $L_{\chi_j}=G$ for $j=18,19$. Choosing identical extensions on the two sides,
the conclusion holds. We complete the proof.
\end{proof}

Summarizing this Section, we have proven the following result.
\begin{theorem}\label{thm: all p}
Let $G={}^2F_4(2)'$ and $p$ be odd. Then $G$ satisfies the inductive NAW condition at the prime $p$.
\end{theorem}
\begin{proof}
For $G={}^2F_4(2)'$. Then its order has the odd prime factor  $p=3,5,13$. For $p=13$, we apply Lemma \ref{lem:3.2} because it has a unique cyclic Sylow $13$-subgroup. For $p=3$, Propositions \ref{p3:prop:equivariant} and \ref{thm:p3} show that Proposition \ref{prop:criterion} holds.  Propositions \ref{p5:prop:equivariant} and \ref{thm:p5} verify the case of $p=5$. We are done.
\end{proof}

\section{Proof of main theorem}\label{section 6}
In the Section, we gives the proof of Theorem \ref{thm:main}.

\begin{proof}[Proof of Theorem \ref{thm:main}]
Let $G$ be a Suzuki group  ${}^2B_2(2^{2m+1})$ for $m\geq 1$, or one of the Ree groups  ${}^2G_2(3^{2m+1})$ and ${}^2F_4(2^{2m+1})$ for $m\geq 0$. \cite[Theorem C]{FFZ} has already addressed  the case of the defining characteristic,
so we assume that $p$ is different from the defining characteristic of $G$. 

For $m\geq1$, by \cite[Proposition 10.9]{FFZ}, Theorems \ref{thm:Suzuki} and \ref{thm:2F_4}, $G$ satisfies the inductive NAW condition at all primes $p$.

For $m=0$, by Theorem \ref{thm: all p} and \cite[Proposition 10.8]{FFZ}, the groups ${}^2F_4(2)'$ and ${}^2G_2(3)'$ satisfy the inductive NAW condition at all primes $p$. We are done.
\end{proof}

\appendix
\section*{Appendix A: Class indexing, power maps and fusion maps}\label{app:maps}
When \(G={}^2F_4(2)'\) and \(L={}^2F_4(2)\). The twenty-two ordinary  class positions of \(G\) in \cite{OrdData} are
\[
\begin{gathered}
(1A,2A,2B,3A,4A,4B,4C,5A,6A,8A,8B,8C,8D,
10A,12A,12B,13A,13B,16A,16B,16C,16D).
\end{gathered}
\]
For \(p=3\), delete positions \(4,9,15,16\); for \(p=5\), delete positions \(8,14\). Thus the $p$-regular class counts are \(18\) and \(20\), respectively.

The third-power and fifth-power maps, given by their images on the full list of positions, are
\begin{align*}
 \pi_3={}&(1,2,3,1,5,6,7,8,3,11,10,12,13,14,5,5,17,18,22,21,20,19),\\
 \pi_5={}&(1,2,3,4,5,6,7,1,9,10,11,12,13,2,16,15,18,17,21,22,19,20).
\end{align*}
These are the third and fifth power-map fields of the ordinary \(G\)-record. The class fusion map of $G$ into the twenty-nine classes of $L$ is
\[
(1,2,3,4,5,6,7,8,9,10,11,12,12,13,14,14,15,15,16,17,16,17).
\]
It is the \texttt{ALF} record from \(G\) to \(L\) in \cite{OrdData} and gives the permutation of outer automorphisms
\[
\alpha=(12\ 13)(15\ 16)(17\ 18)(19\ 21)(20\ 22).
\]
For any Brauer row \(b\), \(b^{\sigma_p}(c)=b(c^p)\). Apply this identity on the $p$-regular classes to \(X_p\), after solving \(D_pX_p=Y_p\), and apply the  outer automorphism permutation of conjugate class in the same way. 

\bigskip\noindent\textbf{\large Acknowledgments.}\
The author would like to thank ChatGPT-6 for some helpful discussions. Specifically, when we discuss the Tits group
\({}^2F_4(2)'\) for $p=5$, we use it to construct local characters in Lemma \ref{p5:lem:localchars}. In addition, we used it to discover the class fusion maps and power maps for $p=3,5$ documented in  Reference \cite{OrdData}.  The mathematical ideas, the proofs, and the writing of the paper were entirely developed by the author.

\phantomsection

\end{document}